\documentclass[preprint]{elsarticle}
\usepackage[margin=1in]{geometry}
\usepackage[sumlimits,intlimits]{mathtools}
\usepackage{amssymb}
\usepackage{color}
\usepackage{caption}
\usepackage{subcaption}
\usepackage{rotating}
\usepackage{float} % allows use of 'H' in position specifier of figures

\let\div\relax %removes definition of \div (the division symbol)
\DeclareMathOperator{\div}{div}

\DeclarePairedDelimiter{\abs}{\lvert}{\rvert}
\DeclarePairedDelimiter{\norm}{\lVert}{\rVert}

\newdefinition{definition}{Definition}
\newtheorem{theorem}{Theorem}
\newdefinition{remark}{Remark}

\begin{document}
%\vspace*{0in}
%\maketitle

\begin{frontmatter}

\title{Uncertain Loading and Quantifying Maximum Energy Concentration within Composite Structures\tnoteref{thanks}}
\tnotetext[thanks]{This work supported in part through NSF Grant DMS-1211066 and by the Air Force Research Laboratory under University of Dayton Research Institute Contract FA8650-10-D-5011.}

\author[lsu,cct]{Robert Lipton\corref{cor2}}
\ead{lipton@math.lsu.edu}

\author[lsu]{Paul Sinz}
\ead{psinz1@lsu.edu}

\author[udri]{Michael Stuebner}
\ead{mstuebner1@udayton.edu}

\cortext[cor2]{Corresponding author}

\address[lsu]{Department of Mathematics, Louisiana State University, Baton Rouge, LA 70803}
\address[cct]{Center for Computation and Technology, Louisiana State University, Baton Rouge, LA 70803}
\address[udri]{University of Dayton Research Institute, Dayton, OH 45469}

\begin{abstract}
We introduce a systematic method for identifying the worst case load among all boundary loads of fixed energy. Here the worst case load  is defined to be the one that delivers the largest fraction of input energy to a prescribed subdomain of interest. The worst case load is identified with the first eigenfunction of a suitably defined eigenvalue problem. The first eigenvalue for this problem is the maximum fraction of boundary energy that can be delivered to the subdomain. We compute worst case boundary loads and associated energy contained inside a prescribed subdomain through the numerical solution of  the eigenvalue problem. We apply this computational method to bound the worst case load associated with an ensemble of random boundary loads given by a second order random process. Several examples are carried out on heterogeneous structures to illustrate the method.
\end{abstract}


\end{frontmatter}

\section{Introduction}
Composite materials often fail near structural features where stress can concentrate. Examples include neighborhoods surrounding lap joints or bolt holes where composite structures are fastened or joined \cite{TongSoutis03}. Large boundary loads deliver energy to the structure and can increase the overall energy near structural features and initiate failure. These considerations provide motivation for a better understanding of energy penetration and concentration inside  structures associated with boundary loading.  One possible approach is to apply the Saint-Venant principle \cite{StVen55}, \cite {Lov27}, \cite{vonMis45}, \cite{MareRuss94} to characterize the rate of decay of the magnitude of the stress or strain away from  the boundary and study its effect on interior subdomains. This type of approach provides  theoretical insight for homogeneous materials. However for composite structures the decay can be slow and far from exponentially decreasing away from the boundary \cite{Horgan}.  With this in mind we attempt a more refined analysis and address the problem from an energy based perspective. In this paper we examine the proportion of the total energy that is contained within a prescribed interior domain of interest in response to boundary displacements or traction loads imposed on the composite structure.

We introduce a computational method for identifying the worst case load defined to be the one that delivers the largest portion of a given input energy to a prescribed interior domain of interest.  The interior domain $\omega$ can surround bolted or bonded joints where stress can concentrate. Here the interior domain  is taken to be a positive distance away from the part of the external boundary of the structural domain  $\omega^*$ where the loads are applied. We show here that it is possible to quantify the  effects of a worst case load that concentrates the greatest proportion of energy onto $\omega$ by a suitably defined {\em concentration} eigenvalue problem.  The largest eigenvalue for the eigenvalue problem  is equal to the maximum fraction of total elastic energy that can be imparted to the subdomain over all boundary loads. The displacement field associated with the worst case load is the eigenfunction associated with the largest eigenvalue. 

As an application we use the concentration eigenvalue problem to bound the fraction of energy imparted on a prescribed subdomain by the worst case load associated with an ensemble of random  loads. While it is possible to consider any type of random boundary loading we illustrate the ideas for boundary loads described by a second order random process  specified in terms of its covariance function and ensemble average.

We conclude noting that related earlier work provides bounds on the local stress and strain amplification generated by material microstructure. Of interest is to identify minimum stress microstructures with the lowest field amplification over all microstructures \cite{SIAPAlaliL},  \cite{ChenLThermoElast}. These results enable the design of graded microstructures for suppression of local stress inside structural components \cite{PhilMagL}, \cite{JMechApplMathLStuebner}.

\section{Energy concentration inside composite structures}\label{sec:problemFormulation}
In this section we develop the notion of energy penetration and its associated concentration within a composite structure. The structural domain $\omega^*$ is  made of a composite material and described by the elastic tensor $\mathbb{A}(x)$ taking different values inside each component material. The composite structures addressed here are general and include fiber reinforced laminates or particle reinforced composites. We suppose that the composite structure is subjected to an ensemble of boundary loads applied to either part or to all of the boundary of the structural domain $\omega^*$. We are interested in the energy concentration around features such as a bolt holes or lap joints contained within a known subdomain $\omega$ of the structural domain $\omega^*$. Here it is assumed that the boundary of  the subdomain $\omega$ is of positive distance away from the part of the structural boundary where loads are being applied.

\begin{figure} [htbp]
   \centering
   \def\svgwidth{2in}
   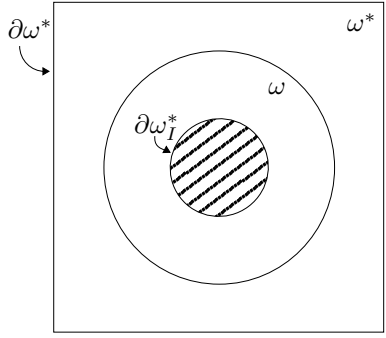
   \caption{Boundary $\partial\omega^*$ and interior domain $\omega$ surrounding bolt hole.}
   \label{Loadandinterior}
\end{figure}

The notion of energy concentration applies to both Dirichlet and traction boundary loading. To fix ideas we first consider Dirichlet loading.  The elastic displacement $u$ is assigned the Dirichlet data $g$ on the exterior boundary  of the domain $\omega^*$ denoted by $\partial\omega^*$. The structural domain $\omega^*$ may be taken to be a bracket or fastener and contain bolt holes  away from the exterior boundary where loads are applied. The boundaries of these holes are assumed clamped and have zero elastic displacement. The collection of these interior boundaries is denoted by $\partial\omega_I^*$, see Figure \ref{Loadandinterior}. The elastic displacement is the solution of the linear elastic system inside the structural domain $\omega^*$ given by
\begin{equation}\label{eq:linearElasticProblem}
\begin{aligned}
   \div\left(\mathbb{A}(x) e(u(x)) \right) = 0,
\end{aligned}
\end{equation}
where $e(u(x))$ is the elastic strain $e(u(x))=(\nabla u(x)+\nabla u(x)^T)/2$ and the elasticity tensor $\mathbb{A}$ satisfies the standard ellipticity and boundedness conditions:
\begin{equation}
\lambda |\zeta|^2\leq \mathbb{A}(x)\zeta:\zeta\leq \Lambda|\zeta|^2,
\label{ellipticity}
\end{equation}
where $\zeta$ is any constant strain tensor, $0<\lambda<\Lambda$ and $\mathbb{A}(x)\zeta:\zeta$ is the elastic energy density given by
\begin{equation}
\mathbb{A}(x)\zeta:\zeta=\sum_{ijkl}\mathbb{A}_{ijkl}(x)\zeta_{kl}\zeta_{ij}
\label{contraction}
\end{equation} 

The work done on the boundary $\partial\omega^*$ delivers the total elastic energy inside the structure and is given by
\begin{equation}
E(g)=\int_{\partial\omega^*}\mathbb{A}e(u)n\cdot g \, ds=\int_{\omega^*} \mathbb{A}e(u):e(u )\, dx,
\label{energy}
\end{equation}
where $n$ is the outward pointing unit normal and $\mathbb{A}e(u)n$ is the traction and $g$ is the boundary displacement.

Now fix a subdomain $\omega$ of interest with boundary a positive distance away from the  structural boundary subjected to loading. Here the domain $\omega$ can be selected to contain a bolt hole or other interior feature of the structure, see Figure \ref{Loadandinterior}.  The elastic energy inside this subdomain is written
\begin{equation}
   (u,u)_{\mathcal{E}(\omega)} = \int_{\omega} \mathbb{A}(x)e(u):e(u) \, dx,
   \label{elasinsubomain}
\end{equation}
and the total elastic energy in the structure is written
\begin{equation}
   (u,u)_{\mathcal{E}(\omega^*)} = \int_{\omega^*} \mathbb{A}e(u):e(u) \, dx.
   \label{elasinstructure}
\end{equation}
The energy norms for the subdomain and structure are written as
\begin{equation}
   \norm{u}_{\mathcal{E}(\omega)}^2 = (u,u)_{\mathcal{E}(\omega)}
   \label{normsubdomain}
\end{equation}
and
\begin{equation}
   \norm{u}_{\mathcal{E}(\omega^*)}^2 = (u,u)_{\mathcal{E}(\omega^*)}=E(g).
   \label{normstructure}
\end{equation}
and the fraction of the total elastic energy in the structure that is contained inside the subdomain $\omega$ is given by the ratio
\begin{equation}
P(g)=\frac{\norm{u}^2_{\mathcal{E}(\omega)}}{E(g)}=\frac{\norm{u}^2_{\mathcal{E}(\omega)}}{\norm{u}^2_{\mathcal{E}(\omega^*)}}.
\label{penetration}
\end{equation}
We define $P(g)$ to be the energy concentration function associated with the boundary displacement $g$.
Collecting our observations we have
\begin{eqnarray}
P(g)=\frac{\norm{u}^2_{\mathcal{E}(\omega)}}{\norm{u}^2_{\mathcal{E}(\omega^*)}}\,:\,\hbox{ with  $\div\left(\mathbb{A}(x) e(u(x)) \right) = 0$ and  $u=g$ on $\partial\omega^*$ }.
\label{penetrationgg}
\end{eqnarray}

We proceed to identify the worst case displacement over all boundary data $g$ associated with  solutions of the elastic system that have square integrable displacements and strain tensors on $\omega^*$. On the interior boundary $\partial\omega_I^*$ we assume clamped boundary conditions for the elastic displacement. We  denote this space of displacements by $H^1_{I0}(\omega^*)$ and the associated  boundary displacements on $\partial\omega^*$ reside in the space  $H^{1/2}(\partial\omega^*)$. It is well known that each boundary displacement in $H^{1/2}(\partial\omega^*)$ corresponds to a unique solution $u$ of the elastic system in $H^1_{I0}(\omega^*)$, see \cite{GiraultRaviart}. Conversely every  solution $u$ of the elastic system in $H^1_{I0}(\omega^*)$ has boundary values in $H^{1/2}(\partial\omega^*)$.
The  worst case displacement $\tilde g$ on $\partial\omega^*$  for the domain $\omega$ delivers the maximum energy concentration
\begin{eqnarray}
P(\tilde g)=\max\left\{P(g)\,:\,\hbox{$g$ in $H^{1/2}(\partial\omega^*)$}\right\}.
\label{worstcasefunctionlargeclass}
\end{eqnarray}

We now state the theorem which identifies the worst case displacement $\tilde g$ in $H^{1/2}(\partial\omega^*)$.

\begin{theorem}
The variational problem given by
\begin{equation}
V=\sup\left\{\frac{\norm{u}^2_{\mathcal{E}(\omega)}}{\norm{u}^2_{\mathcal{E}(\omega^*)}}\,:\,\hbox{$u$ in $H_{I0}^{1}(\omega^*)$ with $\div\left(\mathbb{A}(x) e(u(x)) \right) = 0$}\right\}
\label{energyfracDirichlet}
\end{equation}
has a maximum $\tilde u$ and the maximum energy concentration inside $\omega$ is given by
\begin{eqnarray}
V=P(\tilde g)=\max\left\{P(g)\,:\,\hbox{$g$ in $H^{1/2}(\partial\omega^*)$}\right\}.
\label{worstcasefunction}
\end{eqnarray}
where the worst case displacement is given by the boundary value of $\tilde u$ on $\partial\omega^*$. 
Moreover, the stationary point $\tilde u$ and stationary value $V$ are identified as the eigenfunction eigenvalue pair associated with the solution to the penetration eigenvalue problem given by: 

Find $u$ and $\lambda$ for which $u$ belongs to $H_{I0}^1(\omega^*)$ and is a solution of  $\div\left(\mathbb{A}(x) e(u(x)) \right) = 0$, for $x$ in $\omega^*$, and 
\begin{equation}
   \lambda (u,v)_{\mathcal{E}(\omega^*)} = (u,v)_{\mathcal{E}(\omega)}. 
   \label{eigenproblem}
   \end{equation}
   for  all trial fields $v$ belong to the space of all solutions of $\div\left(\mathbb{A}(x) e(v(x)) \right) = 0$, that belong to $H_{I0}^1(\omega^*)$.
   \label{worstcasetheorem}
\end{theorem}

\begin{remark}
   The maximum $\tilde{u}$ may not be unique.  If the eigenspace corresponding to the largest eigenvalue of \eqref{eigenproblem} has dimension greater than 1, then the maximum $V$ is attained by any vector $u$ in this eigenspace.
\end{remark}

It is clear from the formulation that $V$ does not depend on boundary load but is instead dependent only on the geometry of the boundary of the specimen $\omega^*$, the shape and location of $\omega$, and the microstructure associated with the composite material contained inside the structure. A simple example illustrating the nature of $V$  is given by the  anti-plane shear problem associated with a homogeneous prismatic shaft with circular cross-section $\omega^*$ of radius $1$. The shaft is made from homogeneous isotropic elastic material with specified shear and bulk moduli. On choosing $\omega$ to be a disk of radius $r$ centered inside $\omega^*$ a simple calculation shows that $V=r$. What is notable is that  the maximum fraction of energy that can be delivered to $\omega$ scales as $r$ as opposed to the area fraction of the disk which scales as $r^{2}$.

Next we consider the computation of the worst case traction load and its associated energy concentration. Traction loading $t=\mathbb{A}e(u)n$ is applied to  the boundary  of the domain $\partial\omega^*$ and the  elastic displacement is the solution of the linear elastic system inside the structural domain $\omega^*$ given by
\begin{equation}\label{eq:linearElasticProblemtraction}
\begin{aligned}
   \div\left(\mathbb{A}(x) e(u(x)) \right) = 0,
\end{aligned}
\end{equation}
For this case one considers all traction loads $t$ belonging to the space $H^{-1/2}(\partial\omega^*)$. This space corresponds to all traction boundary values associated with $H^1_{I0}(\omega^*)$ solutions of \eqref{eq:linearElasticProblemtraction}. 
The work done against the applied traction $t$ is equal to the total elastic energy and we write
\begin{equation}
E(t)=\int_{\partial\omega^*} t\cdot u ds=\norm{u}^2_{\mathcal{E}(\omega^*)}.
\label{energytraction}
\end{equation}
The fraction of the total elastic energy in the structure that is contained inside the subdomain $\omega$ is given by the ratio
\begin{equation}
P(t)=\frac{\norm{u}^2_{\mathcal{E}(\omega)}}{E(t)}=\frac{\norm{u}^2_{\mathcal{E}(\omega)}}{\norm{u}^2_{\mathcal{E}(\omega^*)}}.
\label{penetrationt}
\end{equation}
Here we refer to $P(t)$ as the energy concentration function associated with the boundary traction $t$.
The maximum energy concentration inside $\omega$ is given by
\begin{equation}
P(\tilde t)=\max\left\{P(t)\,:\,\hbox{$t$ in $H^{-1/2}(\partial\omega^*)$}\right\}
\label{tractionmaxpenetration}
\end{equation}

We now proceed as before and  state the theorem which identifies the worst case traction $\tilde t$ in $H^{-1/2}(\partial\omega^*)$.

\begin{theorem}
The variational problem given by
\begin{equation}
V=\sup\left\{\frac{\norm{u}^2_{\mathcal{E}(\omega)}}{\norm{u}^2_{\mathcal{E}(\omega^*)}}\,:\,\hbox{$u$ in $H_{I0}^{1}(\omega^*)$ with $\div\left(\mathbb{A}(x) e(u(x)) \right) = 0$}\right\}
\label{energyfracNeumann}
\end{equation}
has a maximum $\tilde u$ and the maximum energy concentration inside $\omega$ is given by
\begin{eqnarray}
V=P(\tilde t)=\max\left\{P(t)\,:\,\hbox{$t$ in $H^{-1/2}(\partial\omega^*)$}\right\}.
\label{worstcasefunctiont}
\end{eqnarray}
and the worst case traction is given by   $\mathbb{A}e(\tilde u)n$ on $\partial\omega^*$. 

Moreover the maximizer $\tilde u$ and $V$ are the same eigenfunction eigenvalue pair associated with the solution to the concentration eigenvalue problem \eqref{eigenproblem}.
\label{worstcasetheoremtraction}
\end{theorem}

\begin{remark}
It is important to note that Theorems \ref{worstcasetheorem} and \ref{worstcasetheoremtraction} apply to the concentration of energy inside $\omega$ associated with work done on the boundary so $\tilde u$ and $V$ are the same for the cases of applied traction and applied displacement.
\end{remark}

In this paper we will compute the worst case boundary loads and associated maximum energy concentration through the numerical solution of  the concentration eigenvalue problem \eqref {eigenproblem}.  Theorems \ref{worstcasetheorem} and \ref{worstcasetheoremtraction}  identify the  space of test and trial functions for use in the computation of the maximum energy concentration $V$ and associated elastic field $\tilde u$. This space is denoted by
$H_\mathbb{A}(\omega^*)$ and is the space of all functions belonging to $H^1_{I0}(\omega^*)$ that are solutions of $\div(\mathbb{A}e(u))=0$.  In Section \ref{sec:numerics1} we provide a numerical method for computing the maximum energy concentration and worst case load by finding $u$ in  $H_\mathbb{A}(\omega^*)$  and $\lambda$ for which
\begin{equation}
   \lambda (u,v)_{\mathcal{E}(\omega^*)} = (u,v)_{\mathcal{E}(\omega)}. 
   \label{eigenproblemnum}
   \end{equation}
   for  all trial fields $v$ belong to $H_\mathbb{A}(\omega^*)$.

\begin{remark}
We conclude noting that the concentration eigenvalue problem \eqref{eigenproblemnum}  corresponds to the largest singular value of the restriction operator $\mathcal{P}$ defined by $\mathcal{P}(u)=u$ for $x$ in $\omega$ where $\mathcal{P}$ acts on the space $H_\mathbb{A}(\omega^*)$, see \cite{Bab11} and \cite{Bab14}. It is shown there that the associated set of singular values decay nearly exponentially and  that the associated eigenfunctions of $\mathcal{P}^*\mathcal{P}$ are a complete orthogonal system for the the space $H_\mathbb{A}(\omega^*)$.
\end{remark}

\section{Energy Concentration due to Random Boundary Data}\label{sec:kkl}

In this section we introduce random loading on the boundary of the domain $\omega^*$  and describe the expected value of the energy concentration associated with the ensemble of random loads. We show that this quantity is bounded above by the  maximum energy concentration obtained through the solution of the eigenvalue problem \eqref{eigenproblemnum}.   It is possible to consider any type of random boundary loading; however, to fix ideas we assume in this treatment that the loading is a second order random process  with specified covariance and expectation. In what follows we apply the Kosambi-Karhunen-Loeve (KKL) expansion (cf. \cite{Gha91}) of the random boundary displacement, $g$ with average, $\bar{g}(x)$, and mean zero fluctuation, $\alpha(x,\theta)$, as
\begin{equation}
   g(x,\theta) = \bar{g}(x) + \alpha(x,\theta).
   \label{randomboundarydata}
\end{equation}
We assume, with no loss of generality, that $\bar{g}(x)$ and $\alpha(x,\theta)$ belong to $H^{1/2}(\partial\omega^*)$.
Where $\alpha(x,\theta)$ is a zero mean, second order stochastic process  with covariance function $C(x_1,x_2)$ defined for  points $x_1$ and $x_2$ on $\partial\omega^*$.  The corresponding KKL expansion of $g(x,\theta)$ is given by
\begin{equation}\label{eq:kklBoundary}
   g(x,\theta) = \bar{g}(x) + \sum_{n=0}^{\infty} \sqrt{\mu_n}\psi_n(\theta)\tilde{g}_n(x)
\end{equation}
where the deterministic functions $\tilde{g}_n(x)$ and parameters $\mu_n$ are the eigenfunctions and eigenvalues of the integral equation
\begin{equation}\label{eq:kklIntegral}
   \int_{\partial\omega^*} C(x_1,x)\hat{g}_n(x,\theta) \, ds = \mu_n \hat{g}_n(x_1,\theta), \hbox{  for $x_1$ on $\partial\omega^*$,}
\end{equation}
where $ds$ is the surface measure with respect to the $x$ variable.
The mean zero random variables, $\psi_n(\theta)$, are determined by 
\begin{equation}
   \psi_n(\theta) = \frac{1}{\sqrt{\mu_n}}\int_{\partial\omega^*} \alpha(x,\theta) \hat{g}_n(x) \, ds
\end{equation}
and the functions $\psi_n(\theta)$ and $\hat{g}_n(x)$ satisfy the orthonormality conditions
\begin{eqnarray}
      &&\langle \psi_m(\theta)\psi_n(\theta) \rangle = \delta_{mn} \label{orthonormalityensb}\\
      &&\int_{\partial\omega^*} \hat{g}_m(x)\cdot\hat{g}_m(x) \, ds = \delta_{mn} \label{orthonormalityenspace}
\end{eqnarray}
for the Kronecker delta function, $\delta_{mn}$, and expectation value, $\langle \cdot \rangle$. The set of all random boundary conditions $g(x,\theta)$ associated with mean $\bar{g}(x)$ and covariance $C(x_1,x_2)$ is denoted by $\mathcal{R}(\partial\omega^*)$.
Using the KKL expansion (\ref{eq:kklBoundary}) in (\ref{eq:linearElasticProblem}), and by linearity, we have an expansion of the solution
\begin{equation}\label{eq:kklSolution}
   u(x,\theta) = \bar{u}(x) + \sum_{n=0}^{\infty} \sqrt{\mu_n}\psi_n(\theta)\hat{u}_n(x)
\end{equation}
where $\bar{u}(x)$ is the solution of 
\begin{equation}\label{eq:kklAverageSolution}
   \left\{
      \begin{aligned}
         \div\mathbb{A}e\bar{u} &= 0  &&\text{ in } \omega^* \\
          \bar{u} &= \bar{g}                   &&\text{ on } \partial \omega^*
      \end{aligned}
   \right.
\end{equation}
and $\hat{u}_n(x,\theta)$ is the solution of
\begin{equation}
   \left\{
      \begin{aligned}
         \div\mathbb{A}e\hat{u}_n &= 0  &&\text{ in } \omega^* \\
         \hat{u}_n &= \hat{g_n}               &&\text{ on } \partial \omega^*
      \end{aligned}
   \right.
\end{equation}
Then $\bar{u}(x)$ is the displacement due to the ensemble average boundary data, $\bar{g}(x)$, and 
\begin{eqnarray}
\sum_{n=0}^{\infty} \sqrt{\mu_n}\psi_n(\theta)\hat{u}_n(x)
\label{fluctuatonsolution}
\end{eqnarray}
 is the displacement due to the random fluctuations in the boundary data.
The boundary data is approximated by truncating (\ref{eq:kklBoundary}) after  $N$ terms giving the approximate solution
\begin{equation}\label{eq:kklApproximateSolution}
   u_N(x,\theta) = \bar{u}(x) + \sum_{n=0}^{N} \sqrt{\mu_n}\psi_n(\theta)\hat{u}_n(x)
\end{equation}
The expectation value of the energy of $u$ over $\omega^*$ is denoted by $\mathcal{E}(\omega^*)$.  For this case we apply the orthonormality condition \eqref{orthonormalityensb} to get
\begin{equation}\label{eq:e*}
   \begin{aligned}
      \mathcal{E}(\omega^*) = \left\langle (u,u)_{\mathcal{E}(\omega^*)} \right\rangle  = \norm{\bar{u}}_{\mathcal{E}(\omega^*)}^2 + \sum_{n=0}^\infty \mu_n \norm{\hat{u}_n}_{\mathcal{E}(\omega^*)}^{2}  
   \end{aligned}
\end{equation}
Truncating the sum after $N$ terms gives the approximation  $u_N$ and the approximate expected value of the energy is given by
\begin{equation}
   \mathcal{E}_N(\omega^*) = \norm{\bar{u}}_{\mathcal{E}(\omega^*)}^2 + \sum_{n=0}^N \mu_n \norm{\hat{u}_n}_{\mathcal{E}(\omega^*)}^{2}  
\end{equation}
The expected value of the energy, $\mathcal{E}(\omega)$, which concentrates within the domain of interest, $\omega$, is defined
\begin{equation}
   \mathcal{E} (\omega)= \left\langle (u,u)_{\mathcal{E}(\omega)} \right\rangle=\norm{\bar{u}}_{\mathcal{E}(\omega)}^2 + \sum_{n=0}^\infty \mu_n \norm{\hat{u}_n}_{\mathcal{E}(\omega)}^{2}  
\end{equation}
with truncation
\begin{equation}
   \mathcal{E}_N(\omega) = \left\langle (u_N,u_N)_{\mathcal{E}(\omega)} \right\rangle=\norm{\bar{u}}_{\mathcal{E}(\omega)}^2 + \sum_{n=0}^N \mu_n \norm{\hat{u}_n}_{\mathcal{E}(\omega)}^{2},
\end{equation}
where the orthogonality conditions \eqref{orthonormalityensb} are used to obtain the right most equalities.
The ratio
\begin{equation}
\overline{P}=\frac{\mathcal{E}(\omega)}{\mathcal{E}(\omega^*)}
\label{averagepenetration}
\end{equation}
is a measure of the expected proportion of energy that is contained within the domain of interest, $\omega$, due to the random boundary conditions. Its truncation $\overline{P}_N$ is defined by
\begin{equation}
\overline{P}_N=\frac{\mathcal{E}_N(\omega)}{\mathcal{E}_N(\omega^*)}.
\label{averagepenetrationN}
\end{equation}

We now consider the maximum energy concentration $P(\mathcal{R})$ into the subdomain $\omega$ associated with the worst possible load ${g}(x,\theta)$ in $\mathcal{R}(\omega^*)$. The maximum energy concentration is defined as
\begin{equation}
P(\mathcal{R})=\max\left\{P(g)\,:\,\hbox{$g=g(x,\theta)$ in $\mathcal{R}(\partial\omega^*)$}\right\}
\label{worstcaserandom}
\end{equation}
where $P(g)$ is the energy concentration associated with the boundary displacement given by \eqref{penetrationgg}.
It is evident from Theorem \ref{worstcasetheorem}   that we have the inequalities
\begin{equation}
\overline{P}\leq P(\mathcal{R})\leq V.
\label{worstcaserandominequalities}
\end{equation}
From the definition of $P(\mathcal{R})$ we have 
\begin{equation}
   \norm{u}^2_{\mathcal{E}(\omega)} \leq P(\mathcal{R}) \norm{u}^2_{\mathcal{E}(\omega^*)}.
\end{equation}
Applying the expectation $\langle \cdot , \cdot \rangle$ the first inequality in \eqref{worstcaserandominequalities} follows.

In the following sections we obtain bounds on the maximum energy concentration $P(\mathcal{R})$  by computing $V$. The maximum energy concentration $V$ is computed through the solution of the eigenvalue problem \eqref{eigenproblemnum}. In Section \ref{sec:numerics2} we compare computed values of the expected energy concentration $\overline{P}$ with the computed values of $V$ for different structural components made from reinforced composite materials.

\section{Computational Approach}\label{sec:numerics1}
In this section we outline our method for computing the expected energy concentration $\overline{P}$ associated with random boundary loading as well as the computation of the maximum energy concentration  $V$ through the solution of the eigenvalue problem \eqref{eigenproblemnum}. We illustrate the method for the antiplane shear problem. For this case  $u$ is the out of plane displacement inside a long prismatic shaft and $\omega^*$ is the shaft cross-section. For the purpose of computation we consider polygonal domains $\omega^*$. The out of plane deformation $u$ is the solution of 
\begin{equation}\label{eq:antiplaneShear}
\left\{
   \begin{aligned}
      \div\left( c(x) \nabla u(x,\theta) \right) = 0 \text{ in } \omega^* \\
      u(x,\theta) = g(x,\theta) \text{ on } \partial \omega^* \\
      u(x,\theta) = 0 \text{ on } \partial \omega_I^*
   \end{aligned}
\right.
\end{equation}
where the coefficient $c(x)$ is the shear modulus taking different values inside each component material.
From (\ref{eq:kklBoundary}) the random boundary displacement is given by $u(x,\theta)=\bar{g}(x)+g(x,\theta)$ with
\begin{equation}\label{eq:kklMarkovian}
   g(x,\theta) = \sum_{n=0}^\infty \sqrt{\mu_n} \psi_n(\theta) \hat{g}_n(x)
\end{equation}
and we choose the Markovian covariance function
\begin{equation}
   C(x_1,x_2) = e^{-\abs{x_1-x_2}/b}
\end{equation}
where $b$ is the correlation length with the same units as $x$.
The domain of interest $\omega$ has boundary a positive distance away from $\partial\omega^*$.  Closed form solutions for the sequence $\{\hat{g}_n\}_{n=1}^\infty$ appearing in (\ref{eq:kklMarkovian}) are obtained through the solution of
\begin{equation}\label{eq:kklIntegralunit}
   \int_{-a}^{a}C(x,x_1)\hat{g}_n(x,\theta) \, dx = \mu_n \hat{g}_n(x_1,\theta), \hbox{  for $x_1$ on $(-a,a)$,}
\end{equation}
as in \cite[page 31]{Gha91} and given by
\begin{equation}
   \begin{aligned}
      \hat{g}_n(x) &= 
            \frac{\cos(\gamma_n x)}{\sqrt{a + \frac{\sin(\gamma_n )}{2\gamma_n}}},  &\mu_n &= \frac{2 c }{\gamma_n^2 + c^2}\\
      \hat{g}'_n(x) &=
            \frac{\sin(\gamma'_n x)}{\sqrt{a - \frac{\sin( \gamma'_n )}{2\gamma'_n}}},  &\mu'_n &= \frac{2c}{{\gamma'}_n^2 + c^2}
   \end{aligned}
   \label{eq:kklSolutions}
\end{equation}
where $c=1/b$ and  $\gamma_n, \gamma'_n \geq 0$ are solutions to the equations
\begin{equation}
   \left\{
   \begin{aligned}
      c-\gamma \tan(\gamma a) &= 0 \\
      \gamma' + c \tan(\gamma a) &= 0.
   \end{aligned}
   \right.
\end{equation}
The random boundary data on $\partial\omega^*$  is provided by mapping the one dimensional KKL solutions \eqref{eq:kklSolutions}  defined on the interval $(-a,a)$ onto the boundary of $\omega^*$.  The expected energy concentration $\overline{P}$ is given by
\begin{equation}
   \overline{P}=\frac{ \norm{\bar{u}}_{\mathcal{E}(\omega)}^2 + \sum_{n=0}^\infty \mu_n \norm{\hat{u}_n}_{\mathcal{E}(\omega)}^{2}}{\norm{\bar{u}}_{\mathcal{E}(\omega^*)}^2 + \sum_{n=0}^\infty \mu_n \norm{\hat{u}_n}_{\mathcal{E}(\omega^*)}^{2}},
   \label{avgenergypenetration}
   \end{equation}
where $\hat{u}_n$ and $\hat{u}'_n$ are the solutions to (\ref{eq:antiplaneShear}) with mapped boundary data $\hat{g}_n$ and $\hat{g}'_n$, respectively. In the simulations we apply the standard finite element method to calculate the expected energy concentration truncated after $N$ terms given by
\begin{equation}
   \overline{P}_N=\frac{ \norm{\bar{u}}_{\mathcal{E}(\omega)}^2 + \sum_{n=0}^N \mu_n \norm{\hat{u}_n}_{\mathcal{E}(\omega)}^{2}}{\norm{\bar{u}}_{\mathcal{E}(\omega^*)}^2 + \sum_{n=0}^N \mu_n \norm{\hat{u}_n}_{\mathcal{E}(\omega^*)}^{2}}.
   \label{avgenergypenetrationN}
   \end{equation}

Next we outline the strategy for numerical solution of the eigenvalue problem \eqref{eigenproblemnum}.  

\begin{enumerate}
\item
In the first step we build a suitably large finite dimensional subspace of $H_\mathbb{A}(\omega^*)$. An $n$-dimensional basis $\{\psi_1,\psi_2,\ldots,\psi_n\}$ of discreet solutions to
\begin{equation}
\div(c(x)\nabla u(x))=0 \hbox{ for $x$ in $\omega^*$}
\nonumber
\end{equation}
is generated from FE solutions associated with a suitable $n$-dimensional set of linearly independent  boundary data posed over $\partial\omega^*$.  Here the $n$-dimensional basis of boundary data is generated by hat functions supported on $\partial\omega^*$. On the interior boundary $\partial\omega^*_I$ the discrete solutions satisfy the clamped boundary condition $u=0$.

\item
The second step is to employ the $n$-dimensional  basis of discrete solutions as $n$-dimensional  test and trial spaces in the finite dimensional discretization of the eigenvalue problem \eqref{eigenproblemnum}. This delivers the finite dimensional generalized eigenvalue problem
\begin{equation}
\lambda A\overline{x}=B\overline{x},
\nonumber
\end{equation}
where $\overline{x}$ is the coordinate vector associated with the basis $\{\psi_1,\psi_2,\ldots,\psi_n\}$ and the stiffness matrices are given by
\begin{equation}
A_{ij}=(\psi_i,\psi_j)_{\mathcal{E}(\omega^*)}\hbox{     and      }B_{ij}=(\psi_i,\psi_j)_{\mathcal{E}(\omega)}.
\nonumber
\end{equation}
\item
The third step solves the generalized eigenvalue problem using standard procedures. Since the matrix $A$ is symmetric and positive definite we  apply a Cholesky factorization,
$A=U^{T}U$, where $U$ is an upper triangular matrix. Then we write the generalized problem as a standard eigenvalue problem $C\overline{y}=\lambda\overline{y}$ where $C=(U^{T})^{-1}BU^{-1}$ and $\overline{y}=U\overline{x}$. This problem is solved using a divide-and-conquer eigenvalue algorithm. The eigenvalues of the standard problem are the same as for the general problem. The eigenvectors $\overline{x}$ are computed by solving $U\overline{x}=\overline{y}$. For all operations routines from the Intel Math Kernel Library \cite{mkl} are used.
\end{enumerate}
The associated eigenfunction $\tilde{u}_n (x)=\sum_{i=1}^n x_i\psi_i(x)$ and eigenvalue $\lambda=V$ deliver the numerical approximation of the worst case load and the maximum energy concentration inside $\omega$.

\section{Numerical Simulations}\label{sec:numerics2}

We carry out numerical approximations of the expected energy concentration $\bar{P}_N$ and the worst case energy concentration $V$ for \eqref{eq:antiplaneShear} over four different geometries with Dirichlet boundary conditions.  For each geometry we specify the parameters $a$, $b$, and $c(x)$ and follow the procedures laid out in Section~\ref{sec:numerics1}.  The parameter $a$ defines the one dimensional domain of the random boundary loads which are mapped to the boundary of each geometry, so $a$ is half the perimeter of the geometry.  Each geometry is defined by a cross-sectional area with a hole and several heterogeneities.  The geometries are defined below in Sections~\ref{subsec:geo1}-\ref{subsec:geo4}. The geometries and the average ensemble load are displayed in Figures~\ref{fig:geo1Layout}-\ref{fig:geo4Layout}.  The computational results for each geometry are displayed in Figures~\ref{fig:geo1}-\ref{fig:geo4}.  All fields portrayed are normalized by the square root of the work done at the boundary of each geometry. Table~\ref{tab:numerics} shows the computed values of $V$ and $\bar{P}_N$ for each geometry as well as the number of functions from the KKL expansion for $g(x,\theta)$ that are used.  All units are dimensionless.  Mixed quadrilateral and triangle elements with bilinear and linear shape functions, respectively, are used for each computation.

\subsection{Geometry 1}\label{subsec:geo1}

\begin{figure}[htbp]
   \centering
   \begin{subfigure}[hb]{0.3\textwidth}
         \includegraphics[width=\textwidth]{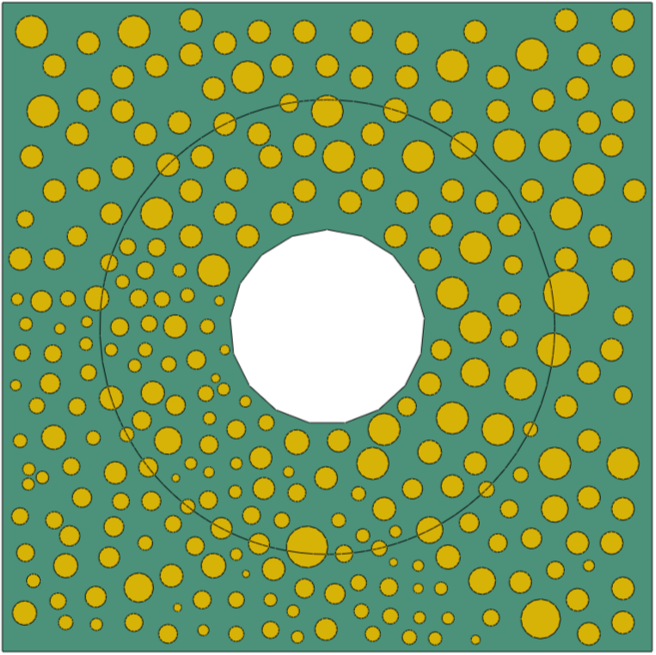}
%      \caption{Cross-sectional layout}
      \caption{}
      \label{fig:geo1Layout}
   \end{subfigure}
\qquad
   \begin{subfigure}[hb]{0.3\textwidth}
      \def\svgwidth{\textwidth}
      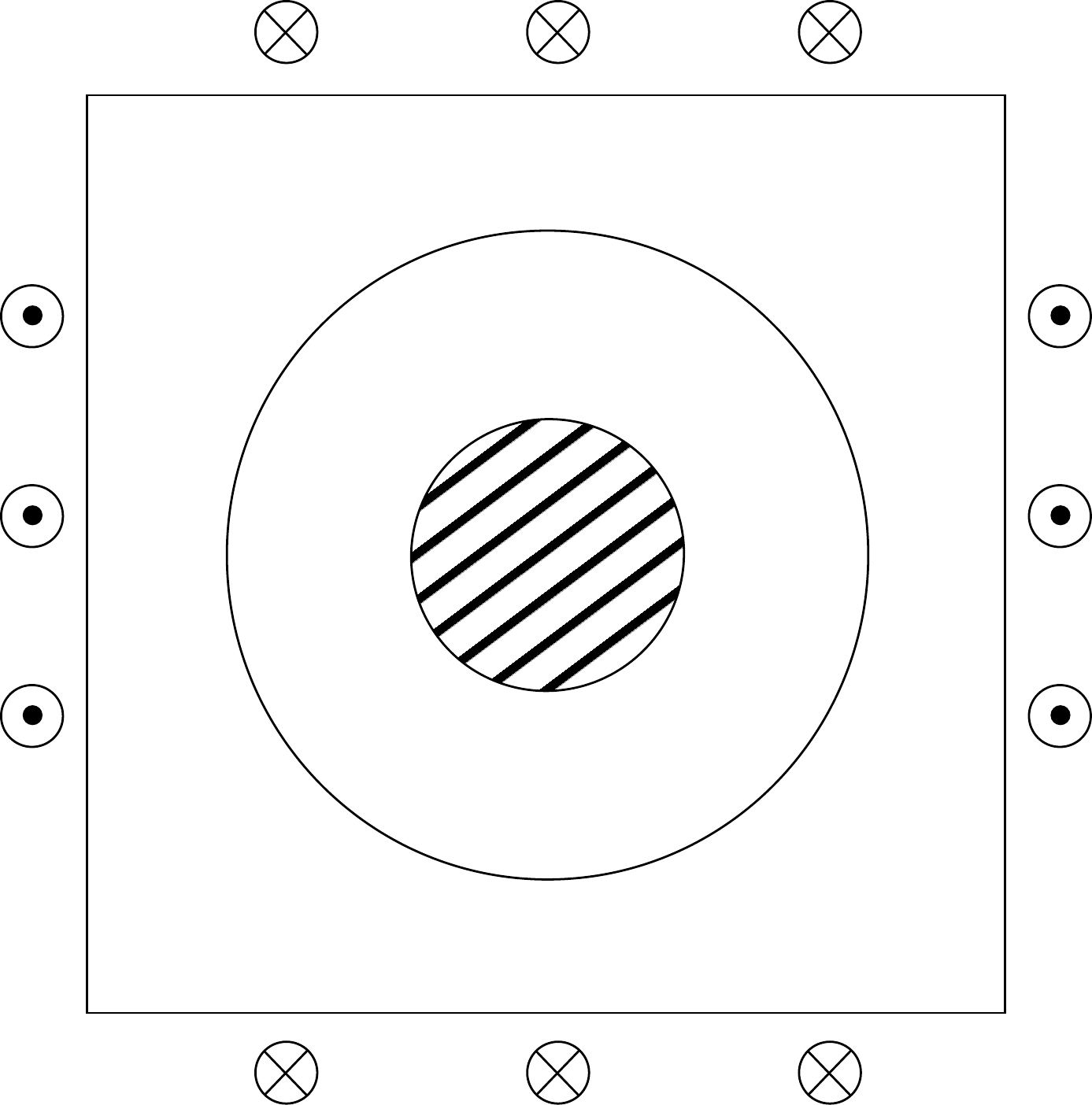
%      \caption{Boundary conditions.  Crosses point into the page, dots point out of the page.}
      \caption{}
      \label{fig:geo1Boundary}
   \end{subfigure}
   \caption{Geometry 1.  (a) The cross sectional layout.  (b) The ensemble averaged Dirichlet boundary conditions for random loading.  Crosses represent displacement into the page, dots out of the page.}
      \label{fig:geo1Definition}
\end{figure}

\begin{figure}[H]
   \centering
   \begin{subfigure}[hb]{0.45\textwidth}
      \includegraphics[trim={4.5in 1.5in 7in 1.5in},clip,width=\textwidth]{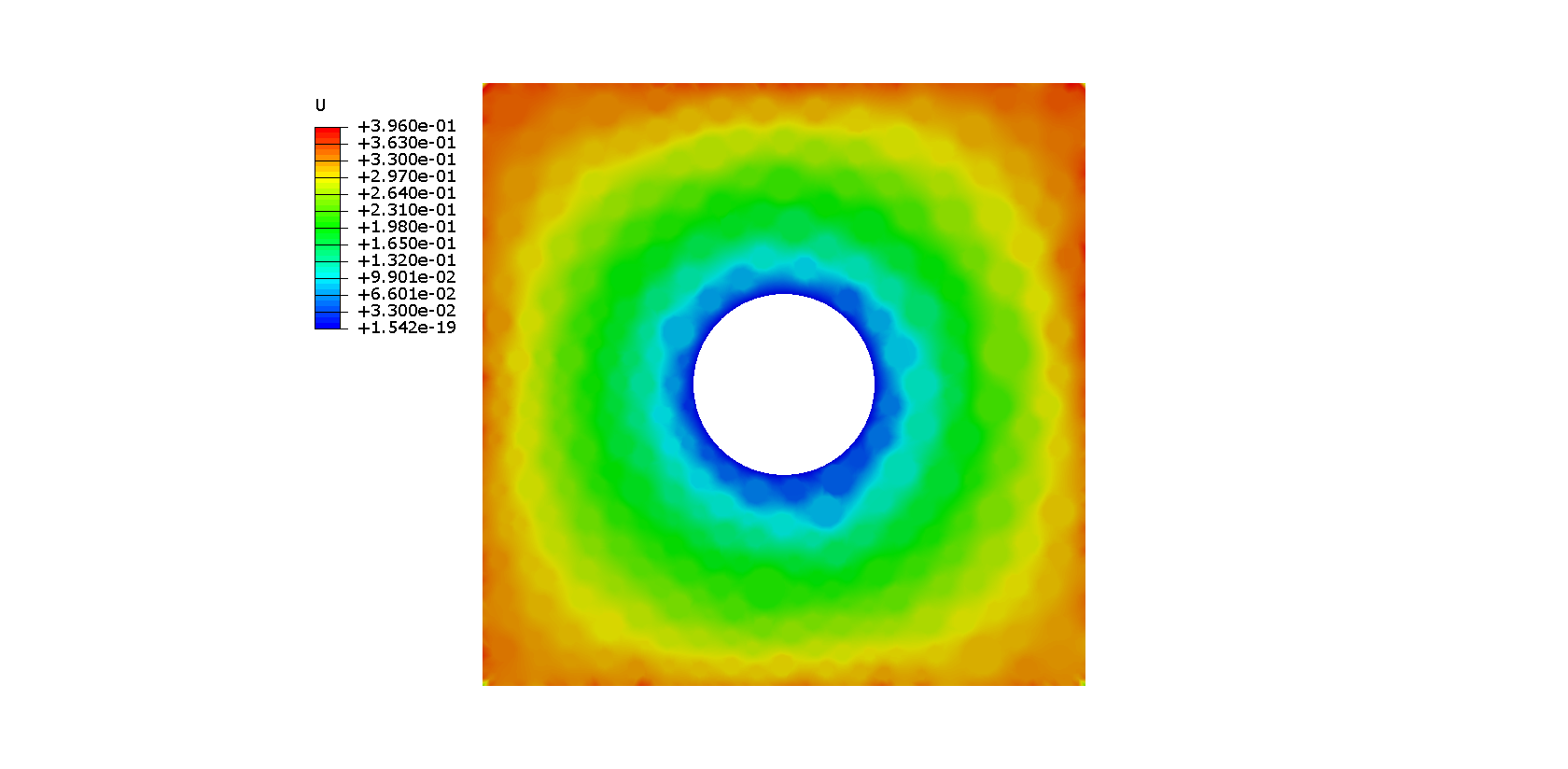}
      \caption{$\tilde{u}(x)$}
      \label{fig:geo1eig}
   \end{subfigure}
   \begin{subfigure}[hb]{0.45\textwidth}
      \includegraphics[trim={4.5in 1.5in 7in 1.5in},clip,width=\textwidth]{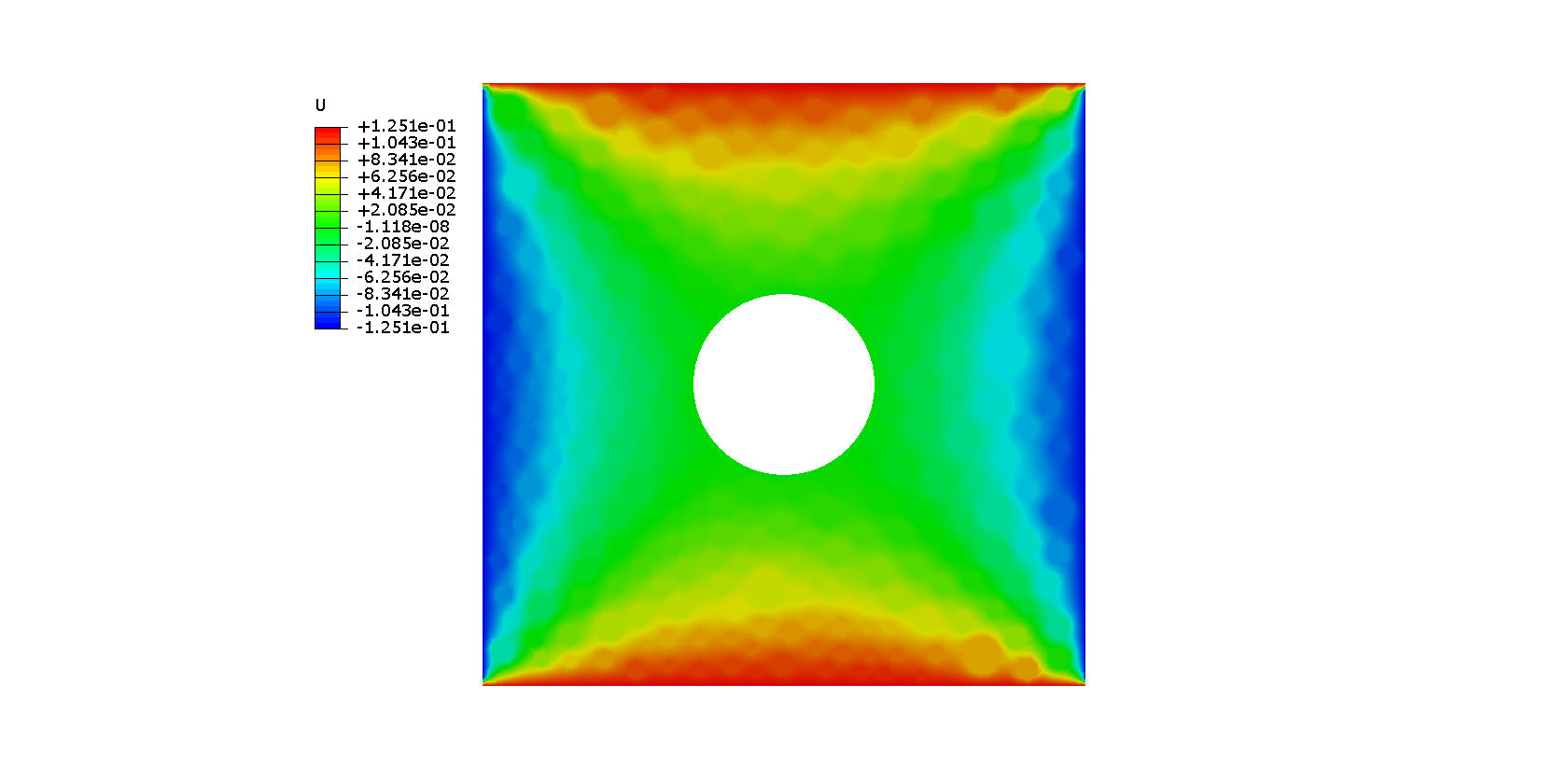}
      \caption{$\bar{u}(x)$}
      \label{fig:geo1kklAvg}
   \end{subfigure}
\\
   \begin{subfigure}[hb]{0.45\textwidth}
      \includegraphics[trim={4.5in 1.5in 7in 1.5in},clip,width=\textwidth]{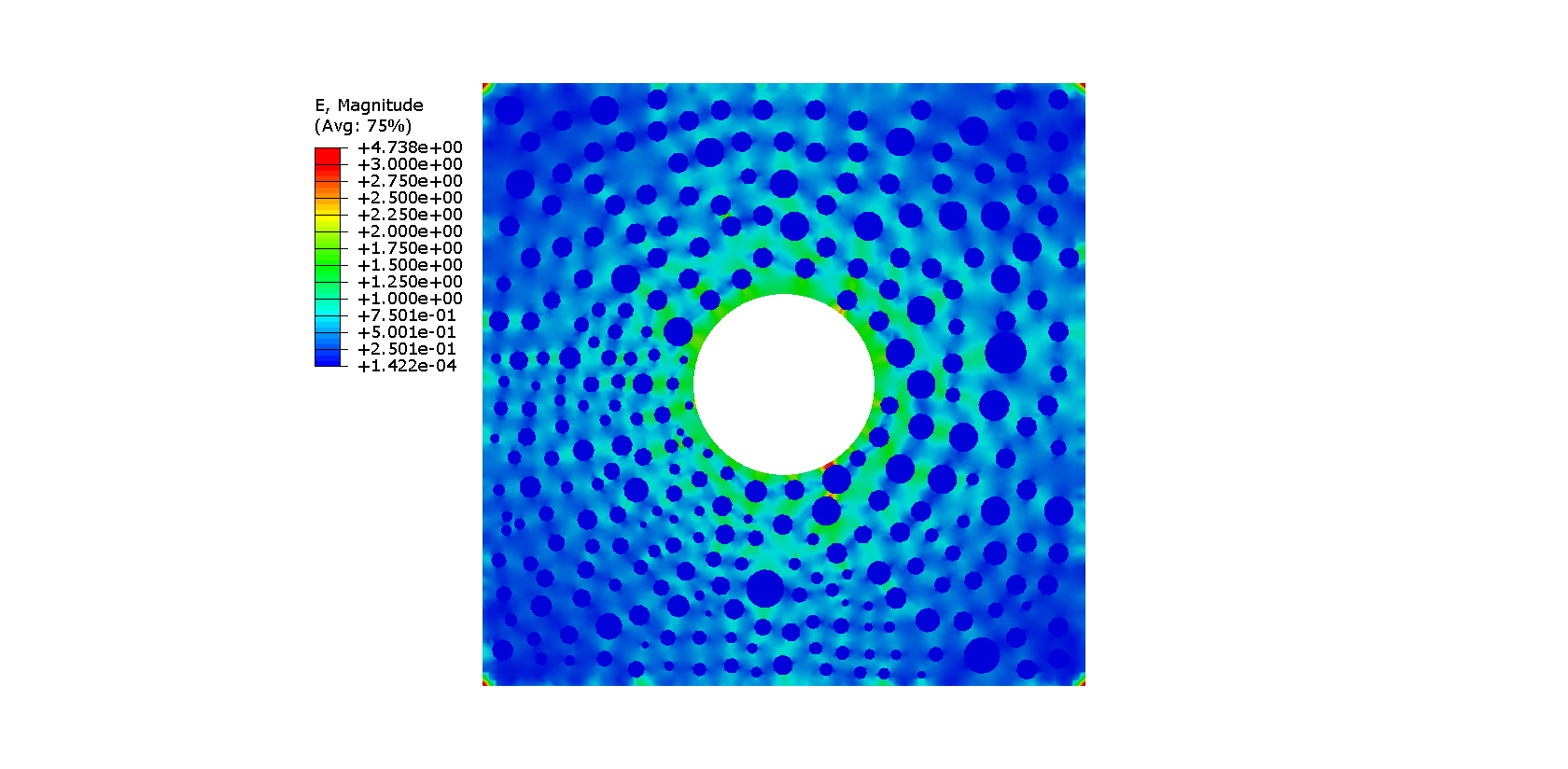}
      \caption{$\abs{\nabla \tilde{u}(x)}$}
      \label{fig:geo1eigStrain}
   \end{subfigure}
   \begin{subfigure}[hb]{0.45\textwidth}
      \includegraphics[trim={4.5in 1.5in 7in 1.5in},clip,width=\textwidth]{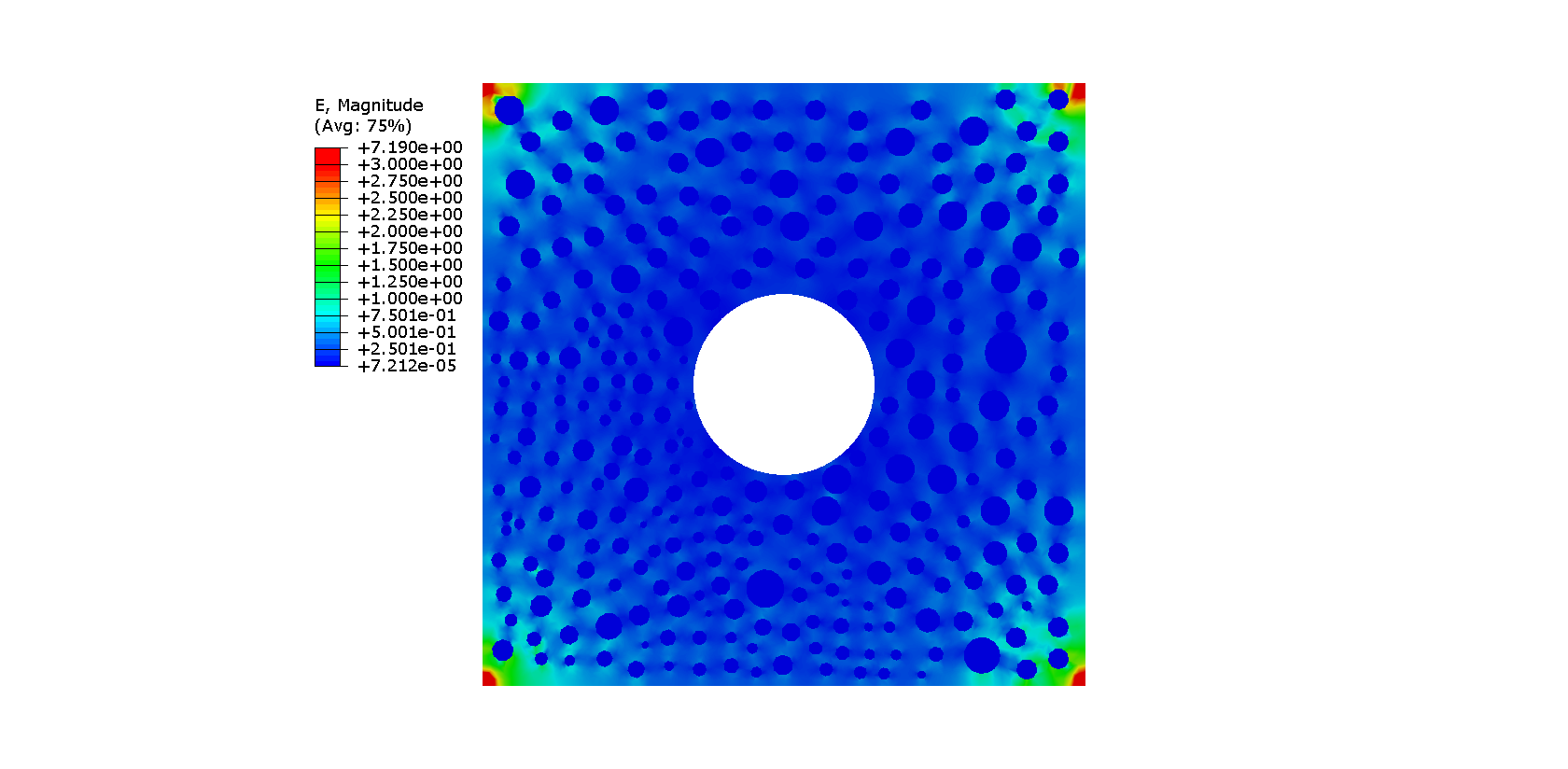}
      \caption{$\abs{\nabla \bar{u}(x)}$}
      \label{fig:geo1kklStrain}
   \end{subfigure}
\\
   \begin{subfigure}[hb]{0.45\textwidth}
      \includegraphics[trim={4.5in 1.5in 7in 1.5in},clip,width=\textwidth]{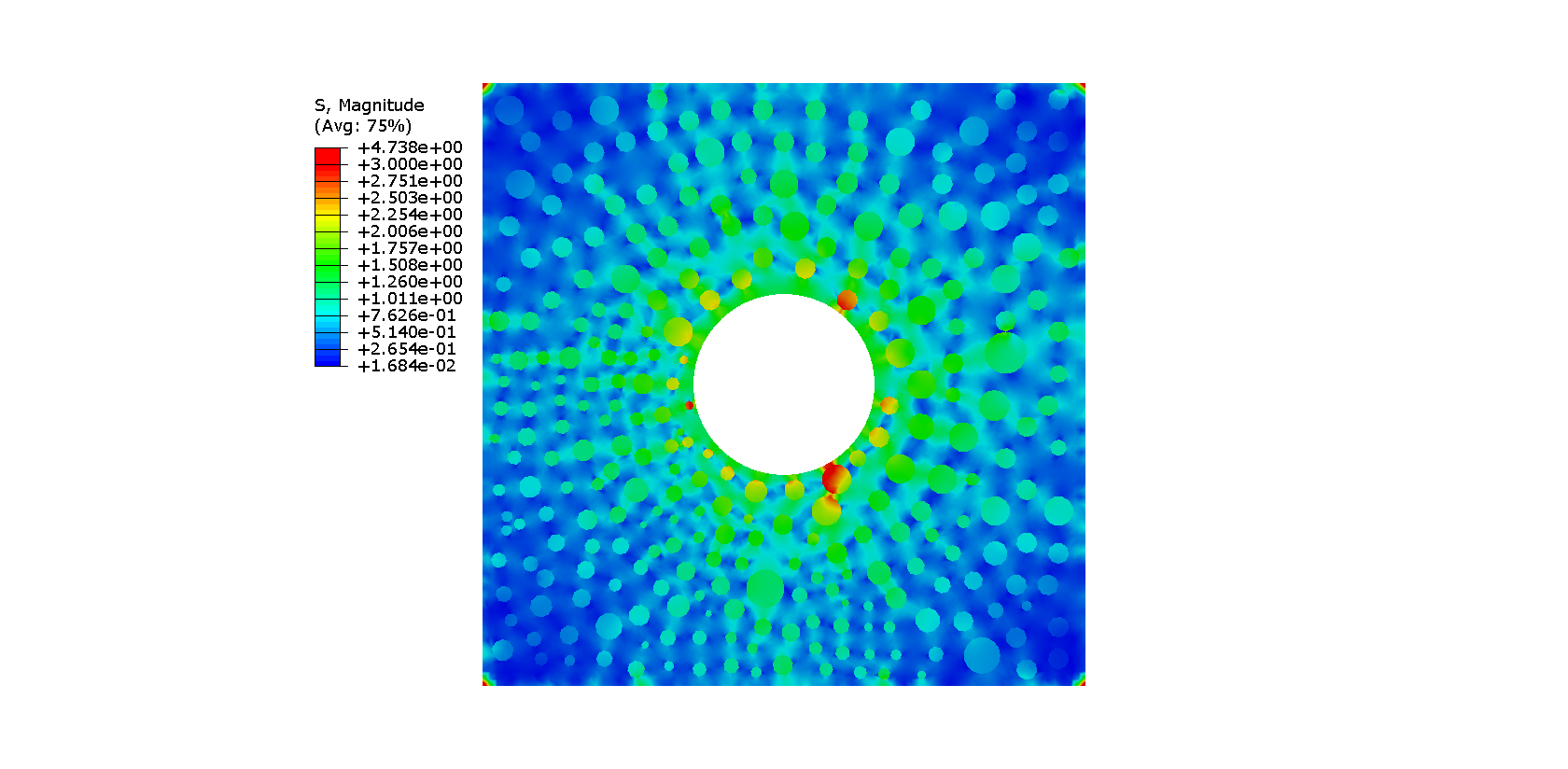}
      \caption{$\abs{c(x)\nabla \tilde{u}(x)}$}
      \label{fig:geo1eigStress}
   \end{subfigure}
   \begin{subfigure}[hb]{0.45\textwidth}
      \includegraphics[trim={4.5in 1.5in 7in 1.5in},clip,width=\textwidth]{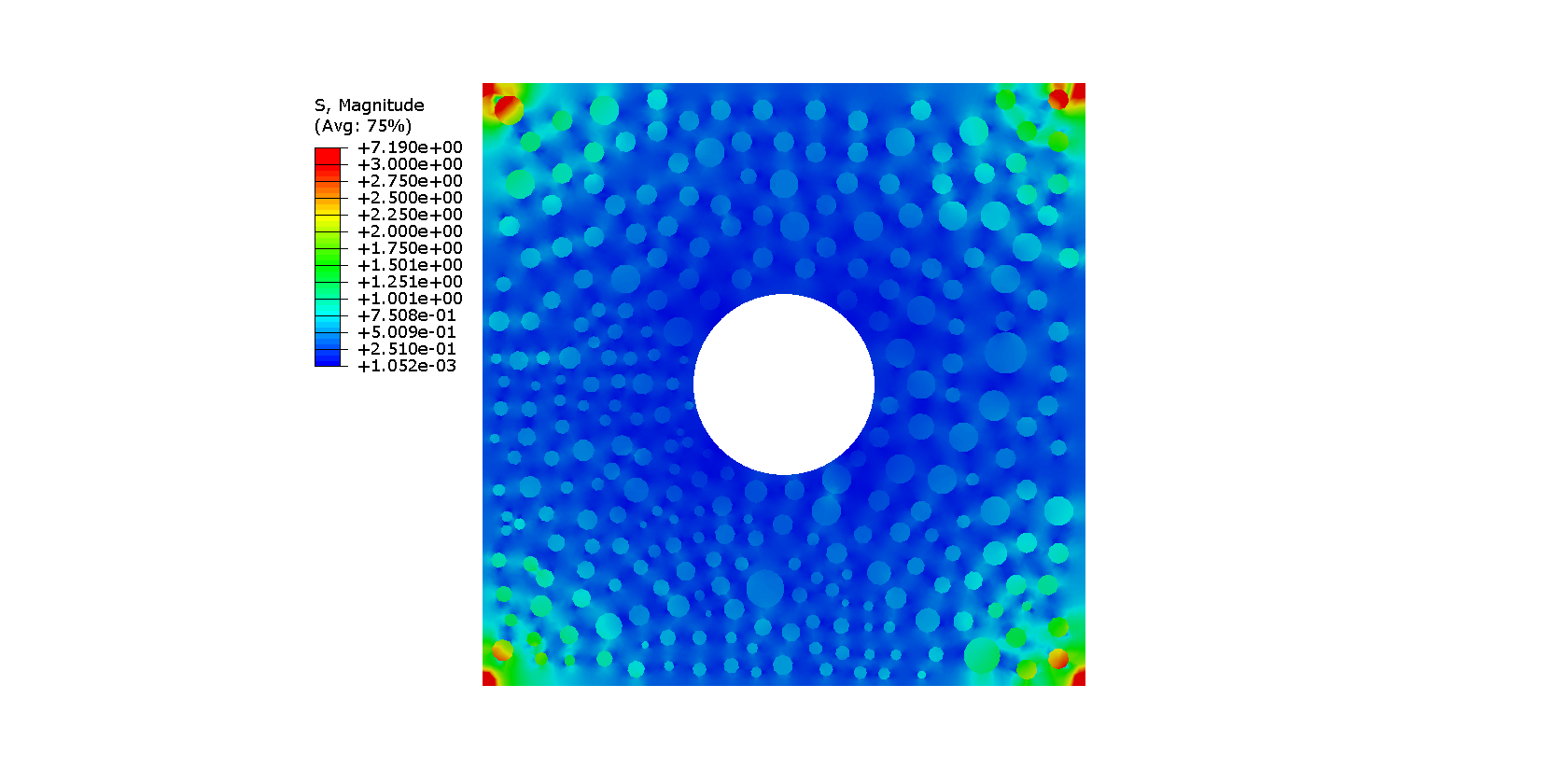}
      \caption{$\abs{c(x)\nabla \bar{u}(x)}$}
      \label{fig:geo1kklStress}
   \end{subfigure}
   \caption{Cross-section of Geometry 1.  The worst case solution (a) $\tilde{u}$ and the corresponding (c) strain and (e) stress compared to the ensemble averaged solution (b) $\bar{u}$ and the corresponding (d) strain and (f) stress.}
   \label{fig:geo1}
\end{figure}

Geometry 1 is a square with a central hole as shown in Figure~\ref{fig:geo1Layout}.  The domain of interest $\omega$ is an annulus around the hole.  The square has side length 2, the hole radius 0.3, and $\omega$ radius 0.7.  The parameter $a$=4, and $b$=1.  The shear modulus $c(x)=1000$ in the inclusions and $c(x)=1$ in the matrix.  The ensemble average boundary displacement $\bar{g}(x)=\pm 0.1$, positive along the edges marked with dotted circles in Figure~\ref{fig:geo1Boundary} and negative along edges marked with crossed circles, and $\bar{g}(x)=0$ along the interior boundary $\partial\omega_I^*$ in Figure~\ref{fig:geo1Boundary}.  In Figure~\ref{fig:geo1} the worst case energy concentration and the associated strain and stress fields are compared to the average ensemble load solution and the associated strain and stress fields. 

\subsection{Geometry 2}\label{subsec:geo2}

\begin{figure}[htbp]
   \centering
   \begin{subfigure}[hb]{0.8in}
      \includegraphics[width=\textwidth]{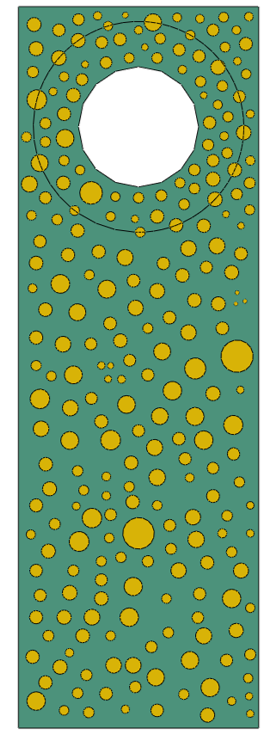}
      \caption{}
      \label{fig:geo2Layout}
   \end{subfigure}
\qquad
   \begin{subfigure}[hb]{0.92in}
      \def\svgwidth{\textwidth}
      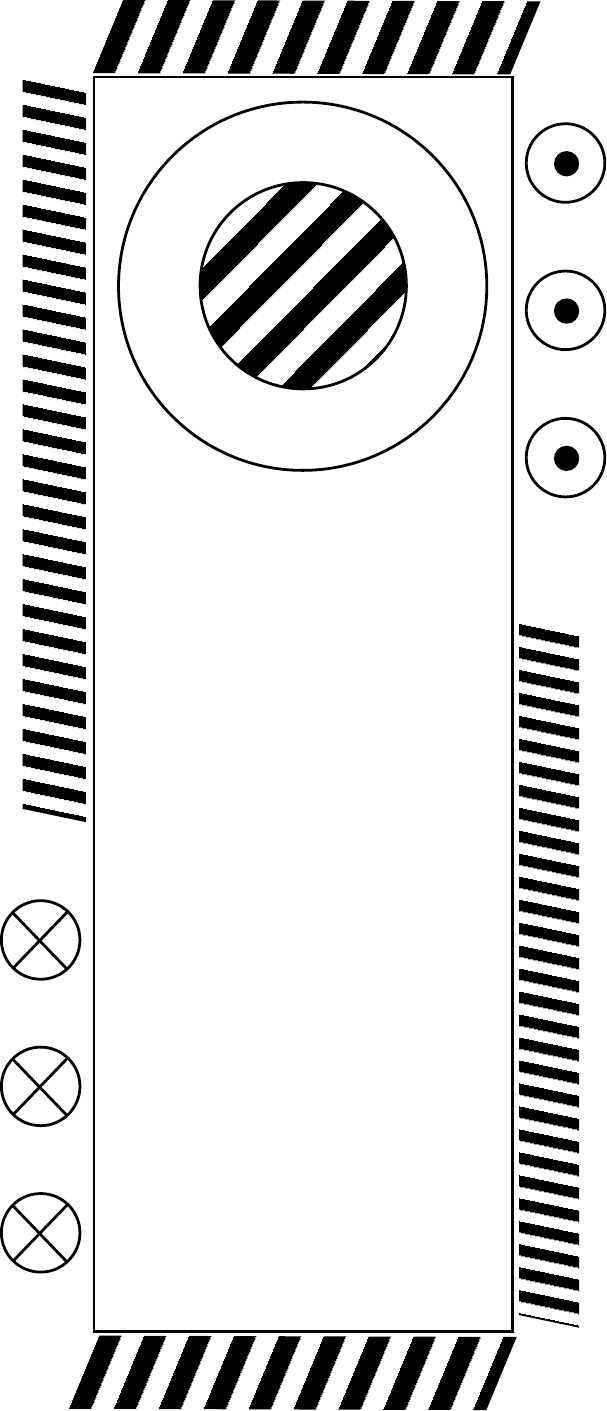
      \caption{}
      \label{fig:geo2Boundary}
   \end{subfigure}
   \caption{Geometry 2.  (a) The cross sectional layout.  (b) The ensemble averaged Dirichlet boundary conditions.  Crosses represent displacement into the page, dots out of the page.}
   \label{fig:geo2Definition}
\end{figure}

Geometry 2 is a rectangle with hole to one end as shown in Figure~\ref{fig:geo2Layout}.  The domain of interest $\omega$ is an annulus around the hole.  The rectangle has height 6 and width 2, the hole radius 0.5, and $\omega$ radius 0.85.  The parameter $a$=8, and $b$=1.  The shear modulus $c(x)=1000$ in the inclusions and $c(x)=1$ in the matrix.  The ensemble average boundary displacement $\bar{g}(x)=\pm 0.1$ or $0$ as shown in Figure~\ref{fig:geo2Boundary}.  Here diagonal lines along edges indicate portions of the outer boundary where $u=0$. The computed fields are displayed in Figure~\ref{fig:geo2}.

\subsection{Geometry 3}\label{subsec:geo3}

\begin{figure}[htpb]
   \centering
   \begin{subfigure}[hb]{0.3\textwidth}
      \includegraphics[width=\textwidth]{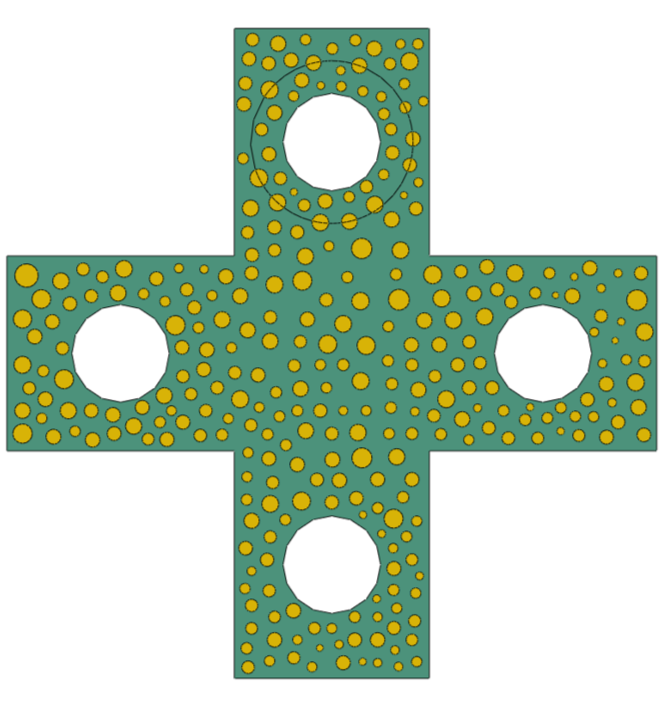}
      \caption{}
      \label{fig:geo3Layout}
   \end{subfigure}
\qquad
   \begin{subfigure}[hb]{0.315\textwidth}
      \def\svgwidth{\textwidth}
      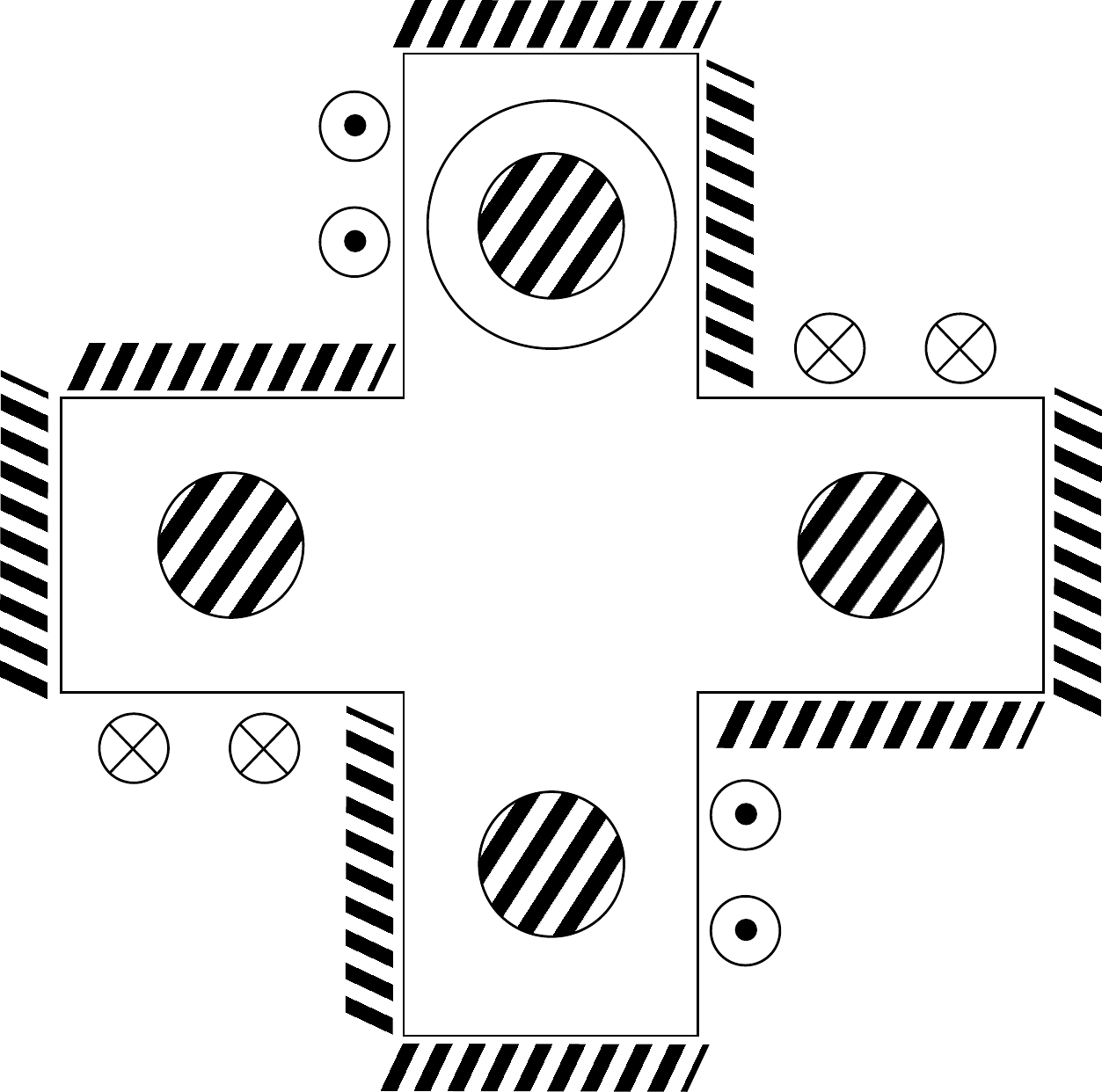
      \caption{}
      \label{fig:geo3Boundary}
   \end{subfigure}
   \caption{Geometry 3.  (a) The cross sectional layout.  (b) The ensemble averaged Dirichlet boundary conditions.  Crosses represent displacement into the page, dots out of the page.}
   \label{fig:geo3Definition}
\end{figure}

\begin{figure}[H]
   \centering
   \begin{subfigure}[hb]{0.37\textwidth}
      \includegraphics[trim={7in 1in 7in 0.4in},clip,width=\textwidth]{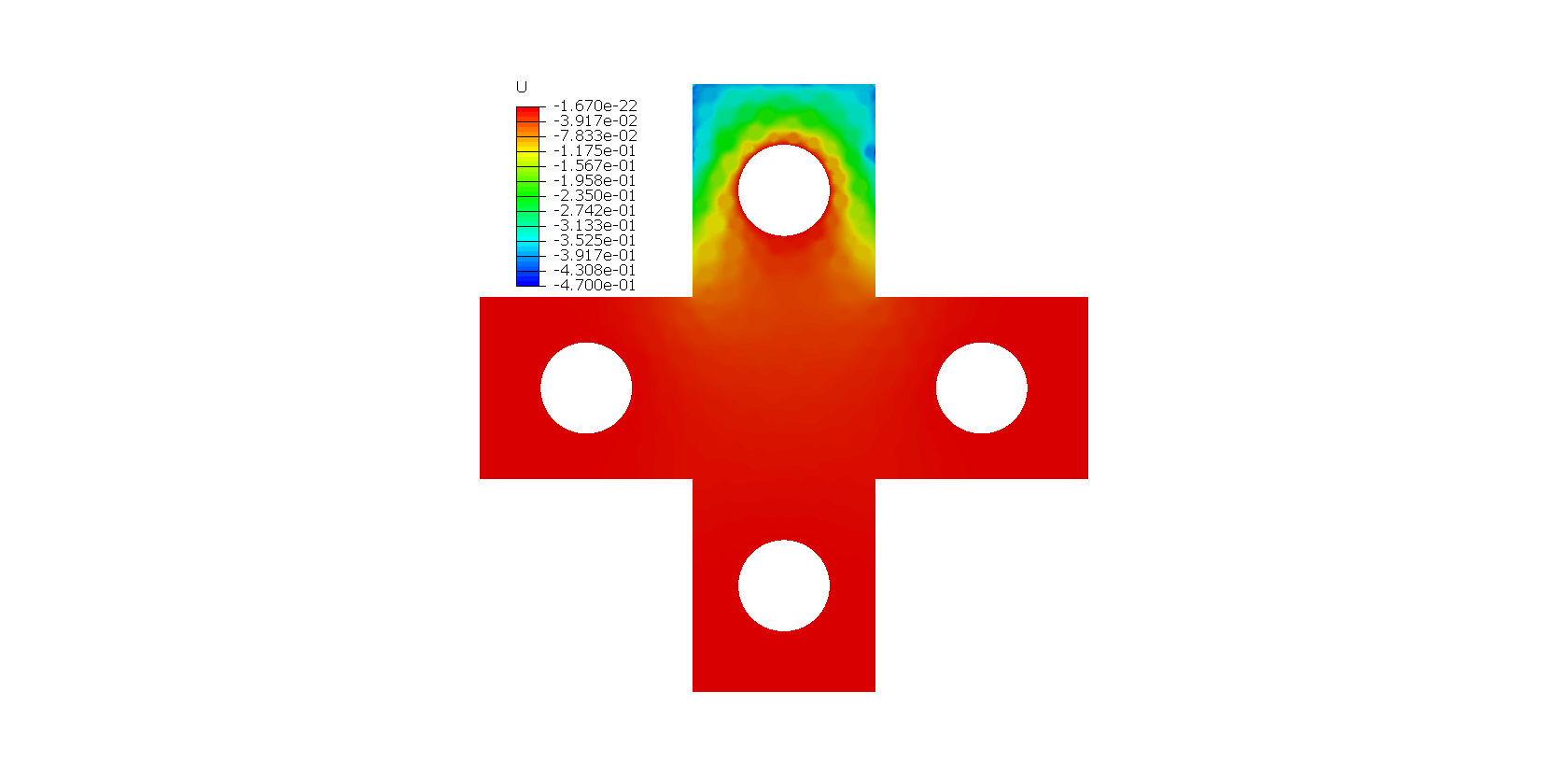}
      \caption{$\tilde{u}(x)$}
      \label{fig:geo3eig}
   \end{subfigure}
   \begin{subfigure}[hb]{0.37\textwidth}
      \includegraphics[trim={7in 1in 7in 0.4in},clip,width=\textwidth]{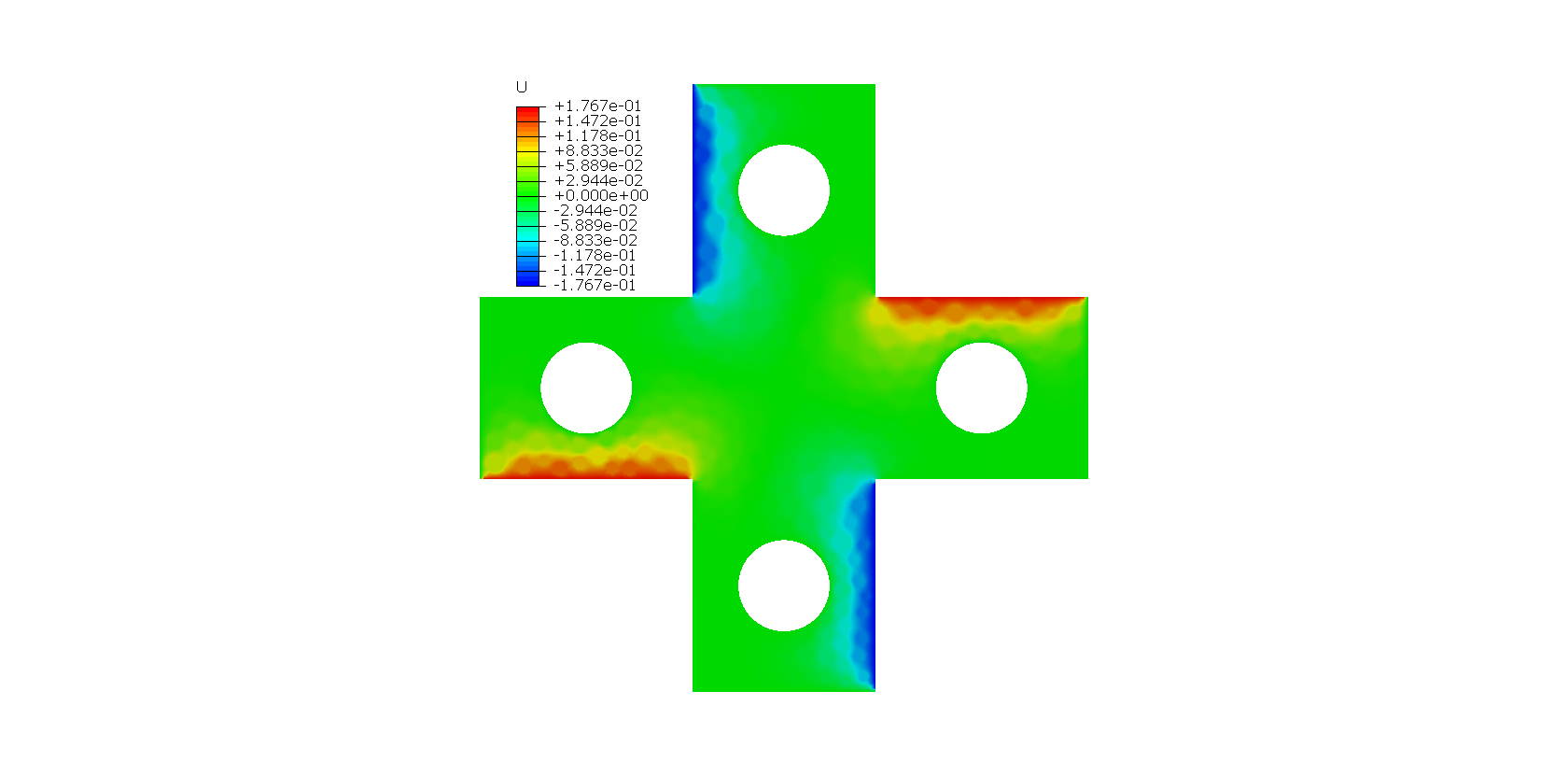}
      \caption{$\bar{u}(x)$}
      \label{fig:geo3kklAvg}
   \end{subfigure}
\\
   \begin{subfigure}[hb]{0.37\textwidth}
      \includegraphics[trim={7in 1in 7in 0.4in},clip,width=\textwidth]{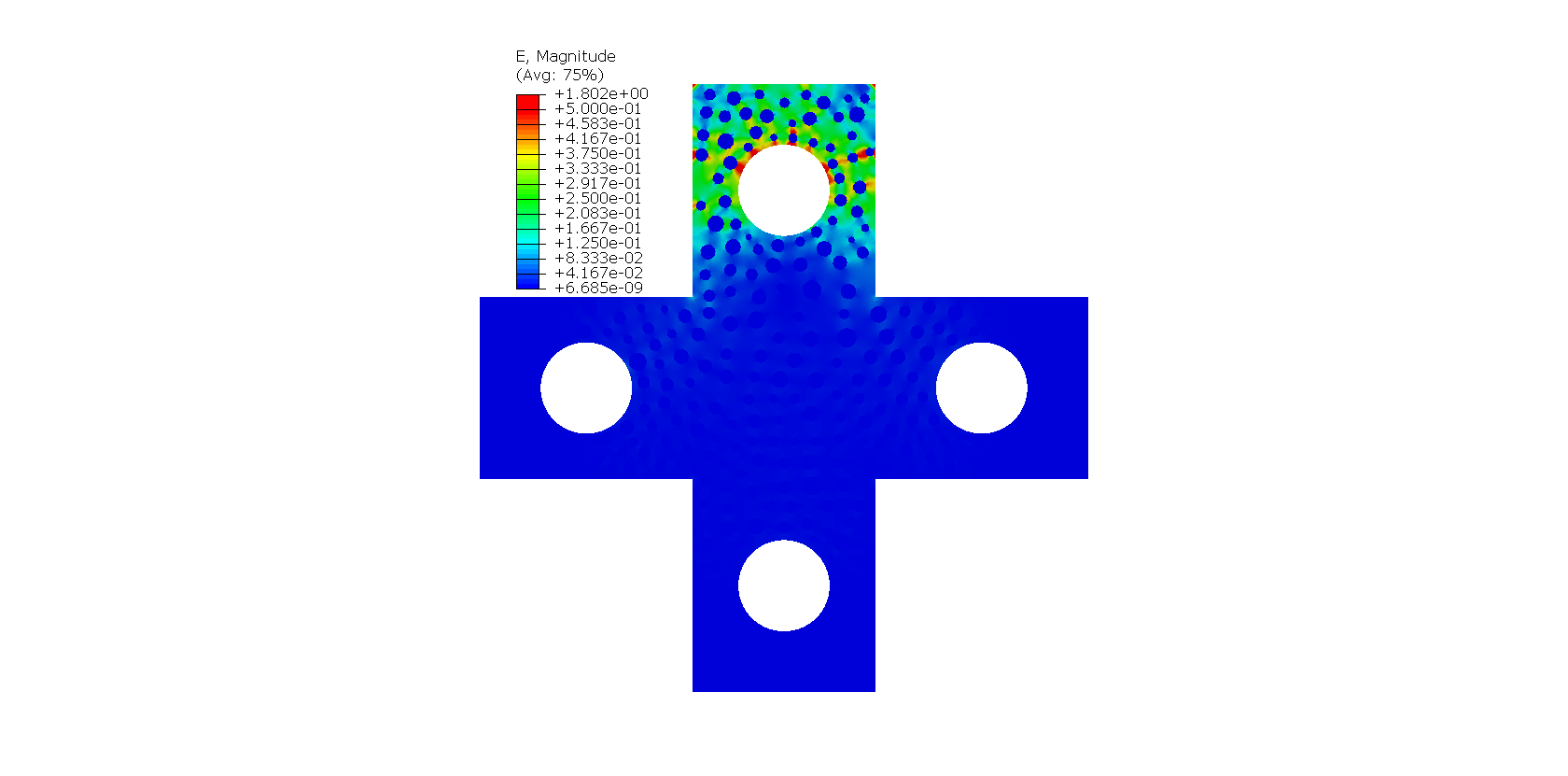}
      \caption{$\abs{\nabla \tilde{u}(x)}$}
      \label{fig:geo3eigStrain}
   \end{subfigure}
   \begin{subfigure}[hb]{0.37\textwidth}
      \includegraphics[trim={5.25in 0.785in 5.25in 0.35in},clip,width=\textwidth]{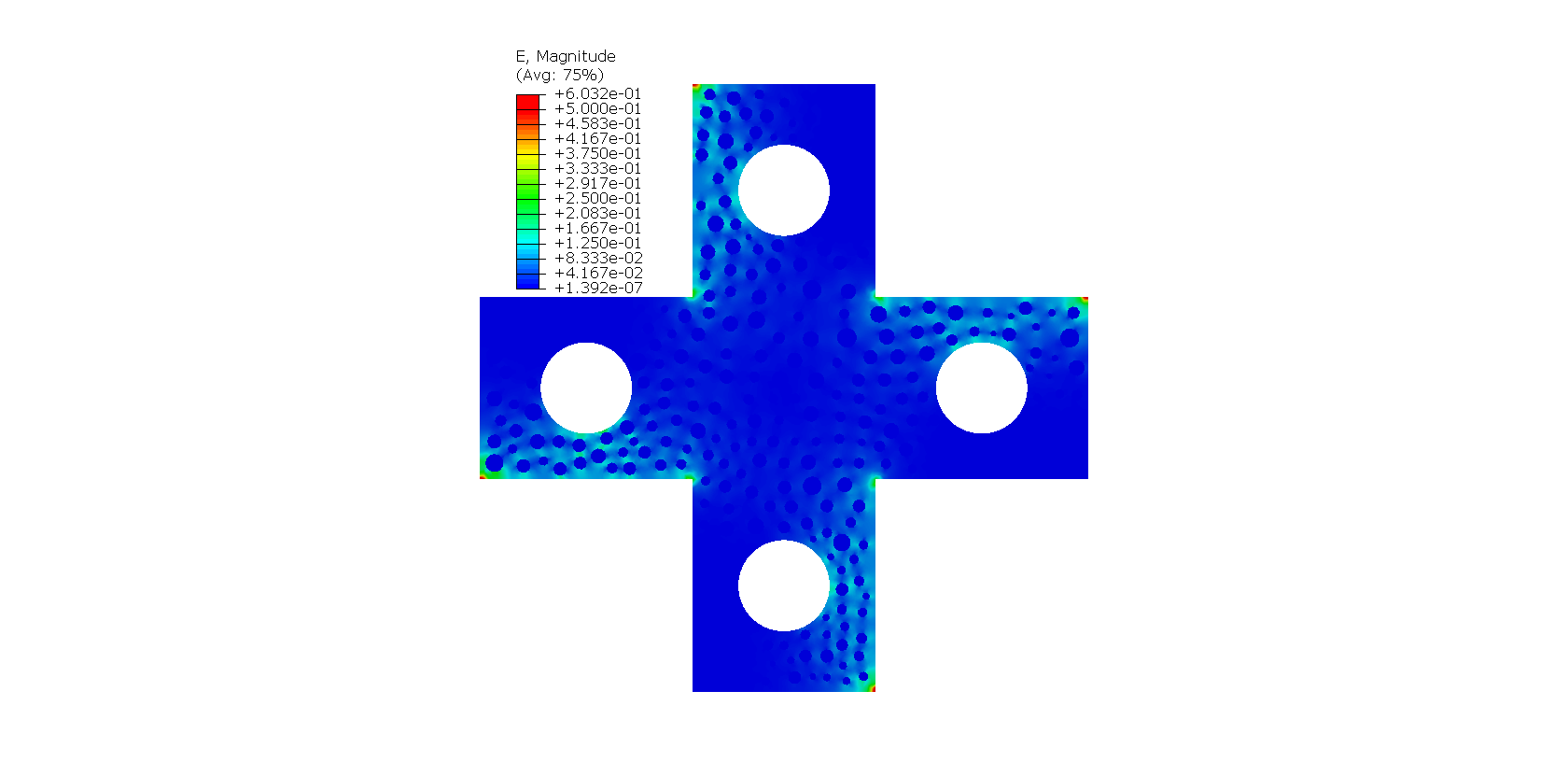}
      \caption{$\abs{\nabla \bar{u}(x)}$}
      \label{fig:geo3kklStrain}
   \end{subfigure}
\\
   \begin{subfigure}[hb]{0.37\textwidth}
      \includegraphics[trim={7in 1in 7in 0.4in},clip,width=\textwidth]{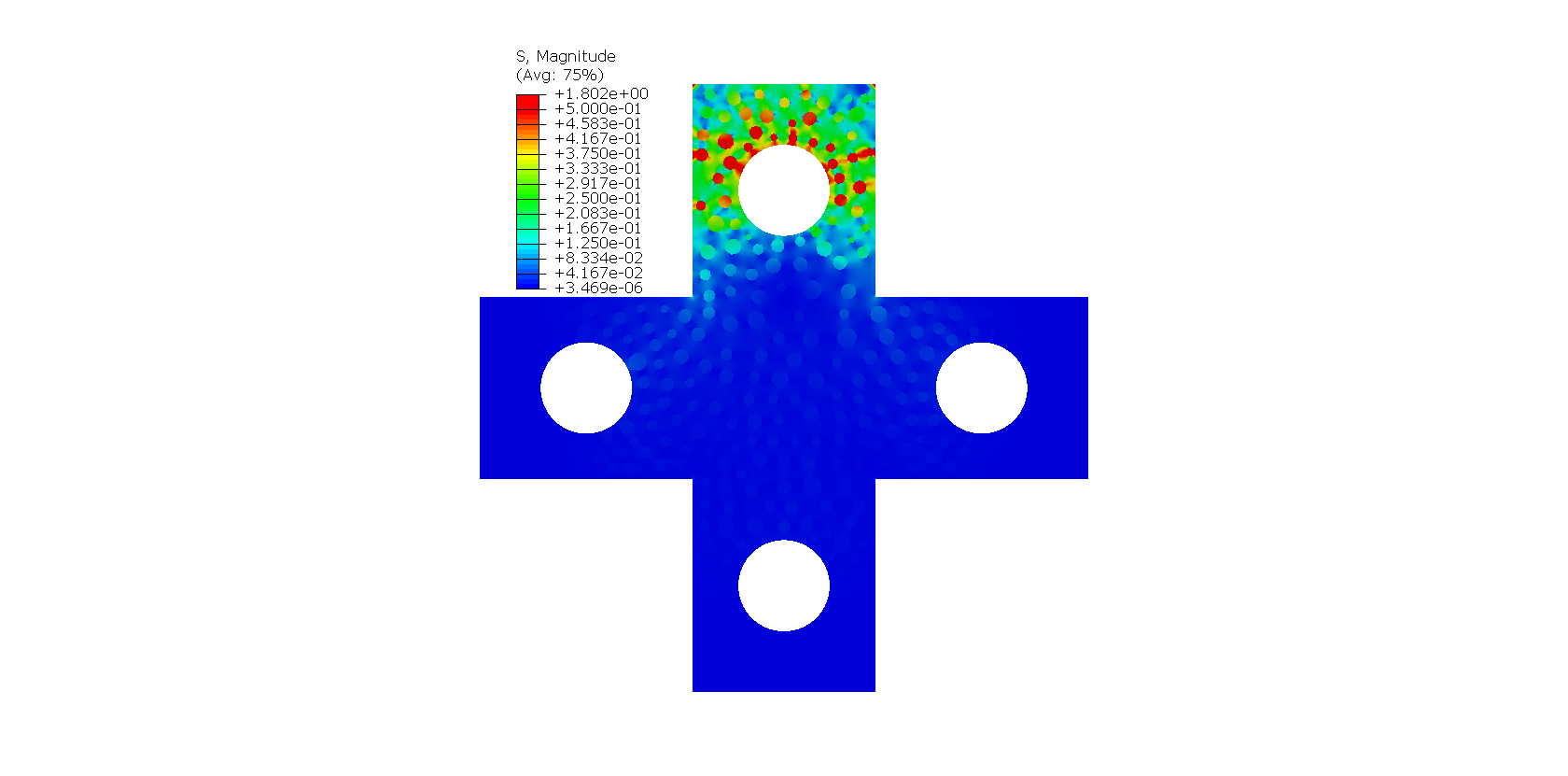}
      \caption{$\abs{c(x)\nabla \tilde{u}(x)}$}
      \label{fig:geo3eigStress}
   \end{subfigure}
   \begin{subfigure}[hb]{0.37\textwidth}
      \includegraphics[trim={7in 1in 7in 0.4in},clip,width=\textwidth]{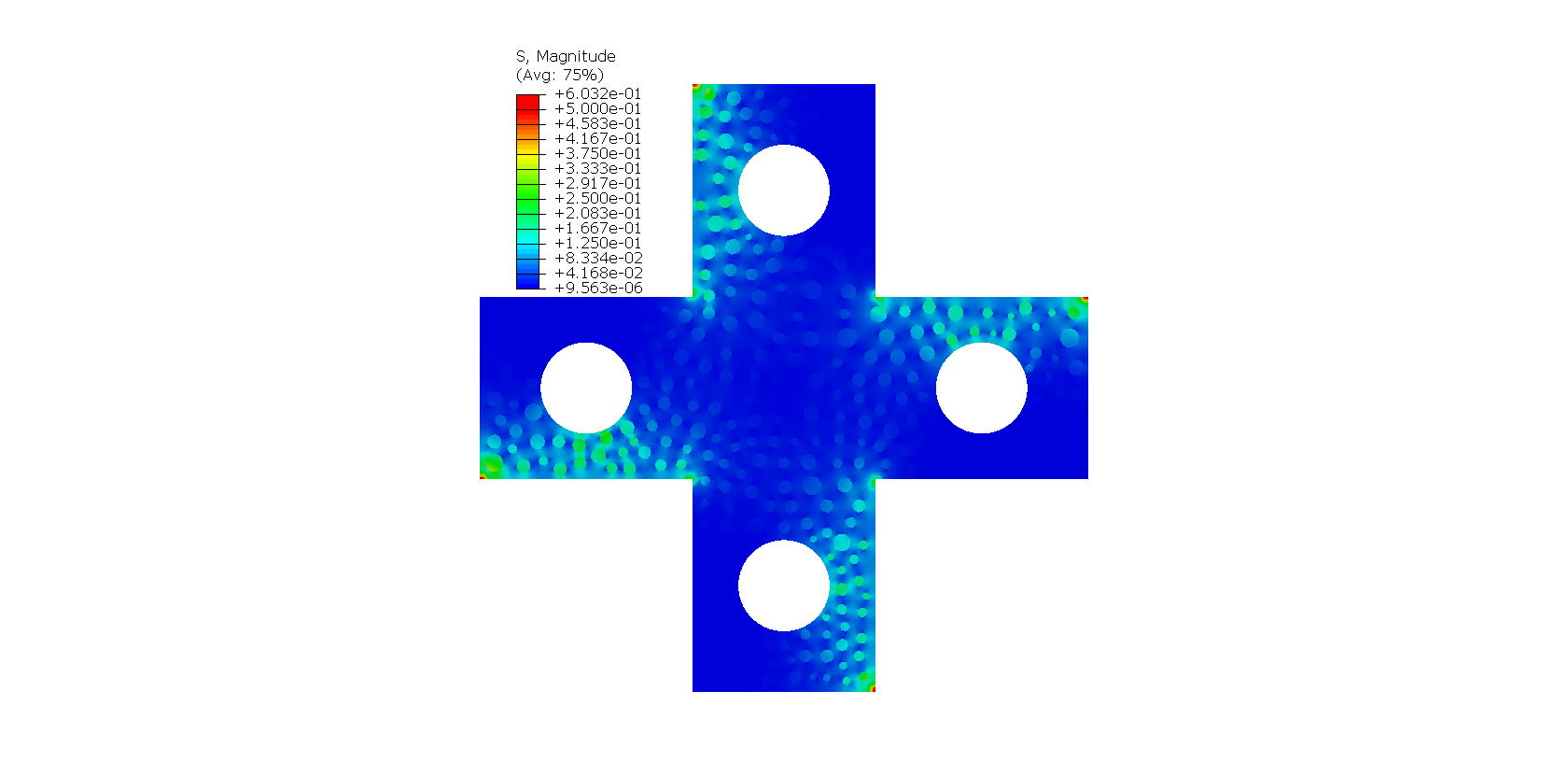}
      \caption{$\abs{c(x)\nabla \bar{u}(x)}$}
      \label{fig:geo3kklStress}
   \end{subfigure}
   \caption{Cross-section of Geometry 3.  The worst case solution (a) $\tilde{u}$ and the corresponding (c) strain and (e) stress compared to the ensemble averaged solution (b) $\bar{u}$ and the corresponding (d) strain and (f) stress.}
   \label{fig:geo3}
\end{figure}

Geometry 3 is a cross with holes at all four ends as shown in Figure~\ref{fig:geo3Layout}.  The domain of interest $\omega$ is an annulus around the top hole.  The cross has a central square of side length 6 and arms of length 7 and width 6, the holes have radius 1.5, and $\omega$ radius 2.5.  The parameter $a$=40, and $b$=1.  The shear modulus $c(x)=1000$ in the inclusions and $c(x)=1$ in the matrix.  The ensemble average boundary displacement $\bar{g}(x)=\pm 0.1$ or $0$ as shown in Figure~\ref{fig:geo3Boundary}.  The computed fields are displayed in Figure~\ref{fig:geo3}.

\subsection{Geometry 4}\label{subsec:geo4}

\begin{figure}[htpb]
   \centering
   \begin{subfigure}[hb]{0.3\textwidth}
      \includegraphics[width=\textwidth]{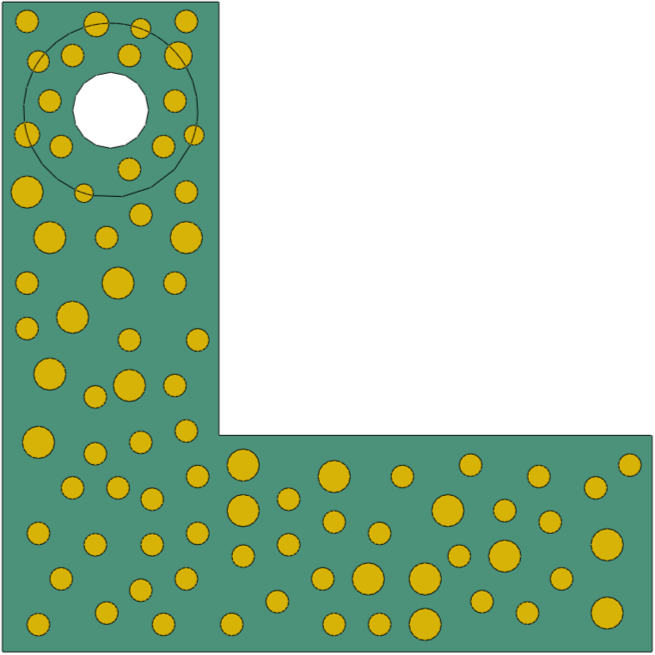}
      \caption{}
      \label{fig:geo4Layout}
   \end{subfigure}
\qquad
   \begin{subfigure}[hb]{0.3\textwidth}
      \def\svgwidth{\textwidth}
      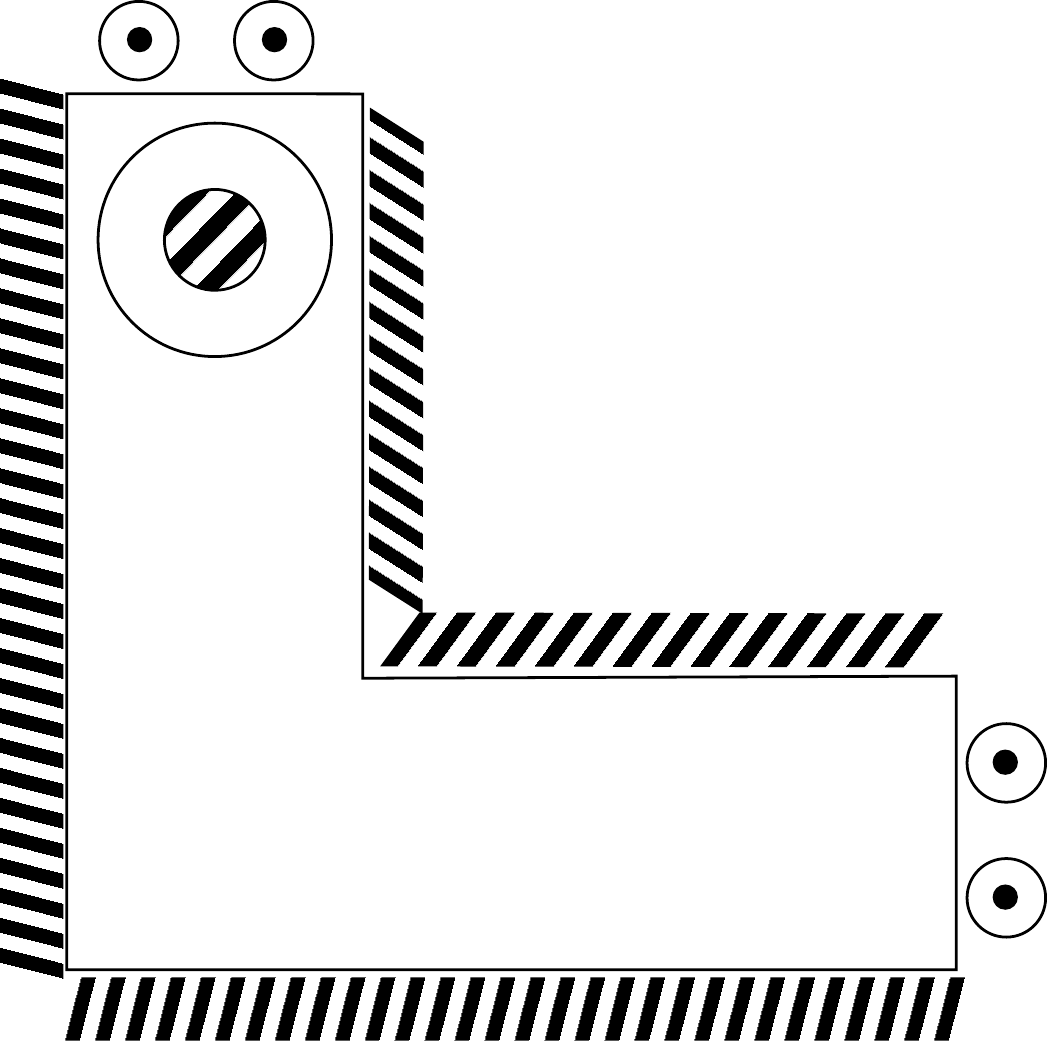
      \caption{}
      \label{fig:geo4Boundary}
   \end{subfigure}
   \caption{Geometry 4.  (a) The cross sectional layout.  (b) The ensemble averaged Dirichlet boundary conditions.  Crosses represent displacement into the page, dots out of the page.}
   \label{fig:geo4Definition}
\end{figure}

\begin{figure}[H]
   \centering
   \begin{subfigure}[hb]{0.39\textwidth}
      \includegraphics[trim={7in 1in 7in 1in},clip,width=\textwidth]{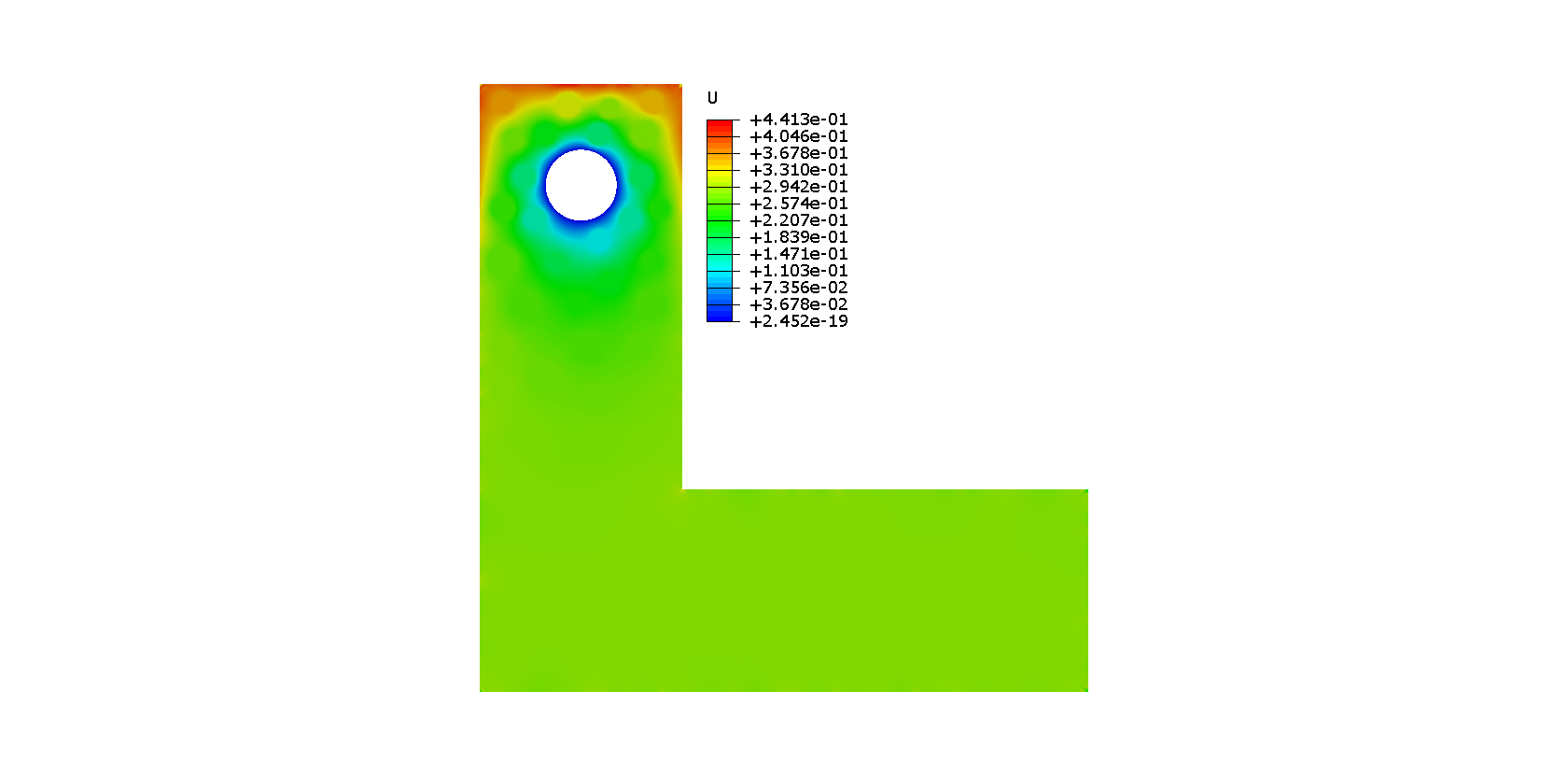}
      \caption{$\tilde{u}(x)$}
      \label{fig:geo4eig}
   \end{subfigure}
   \begin{subfigure}[hb]{0.39\textwidth}
      \includegraphics[trim={7in 1in 7in 1in},clip,width=\textwidth]{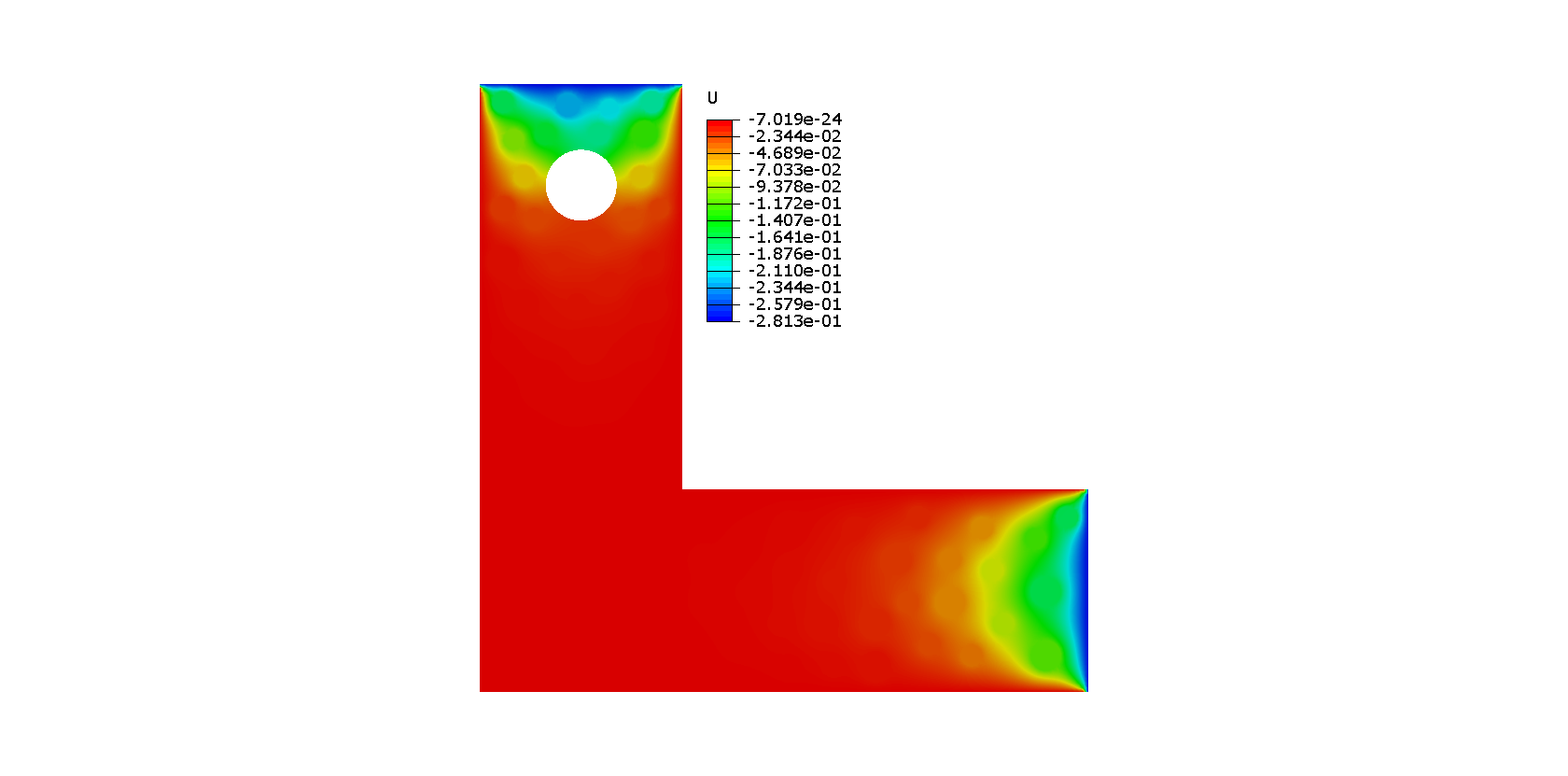}
      \caption{$\bar{u}(x)$}
      \label{fig:geo4kklAvg}
   \end{subfigure}
\\
   \begin{subfigure}[hb]{0.39\textwidth}
      \includegraphics[trim={7in 1in 7in 1in},clip,width=\textwidth]{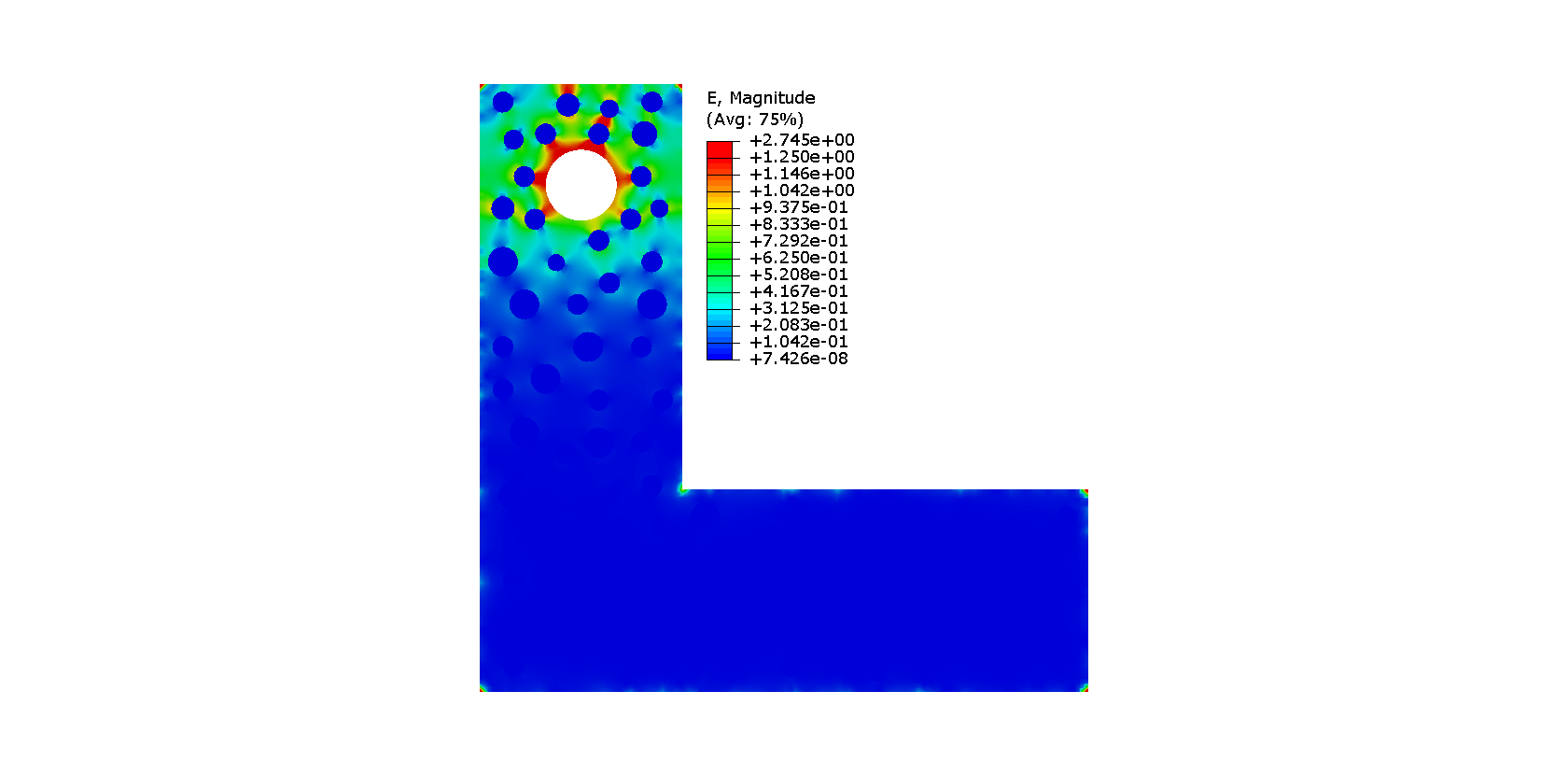}
      \caption{$\abs{\nabla \tilde{u}(x)}$}
      \label{fig:geo4eigStrain}
   \end{subfigure}
   \begin{subfigure}[hb]{0.39\textwidth}
      \includegraphics[trim={7in 1in 7in 1in},clip,width=\textwidth]{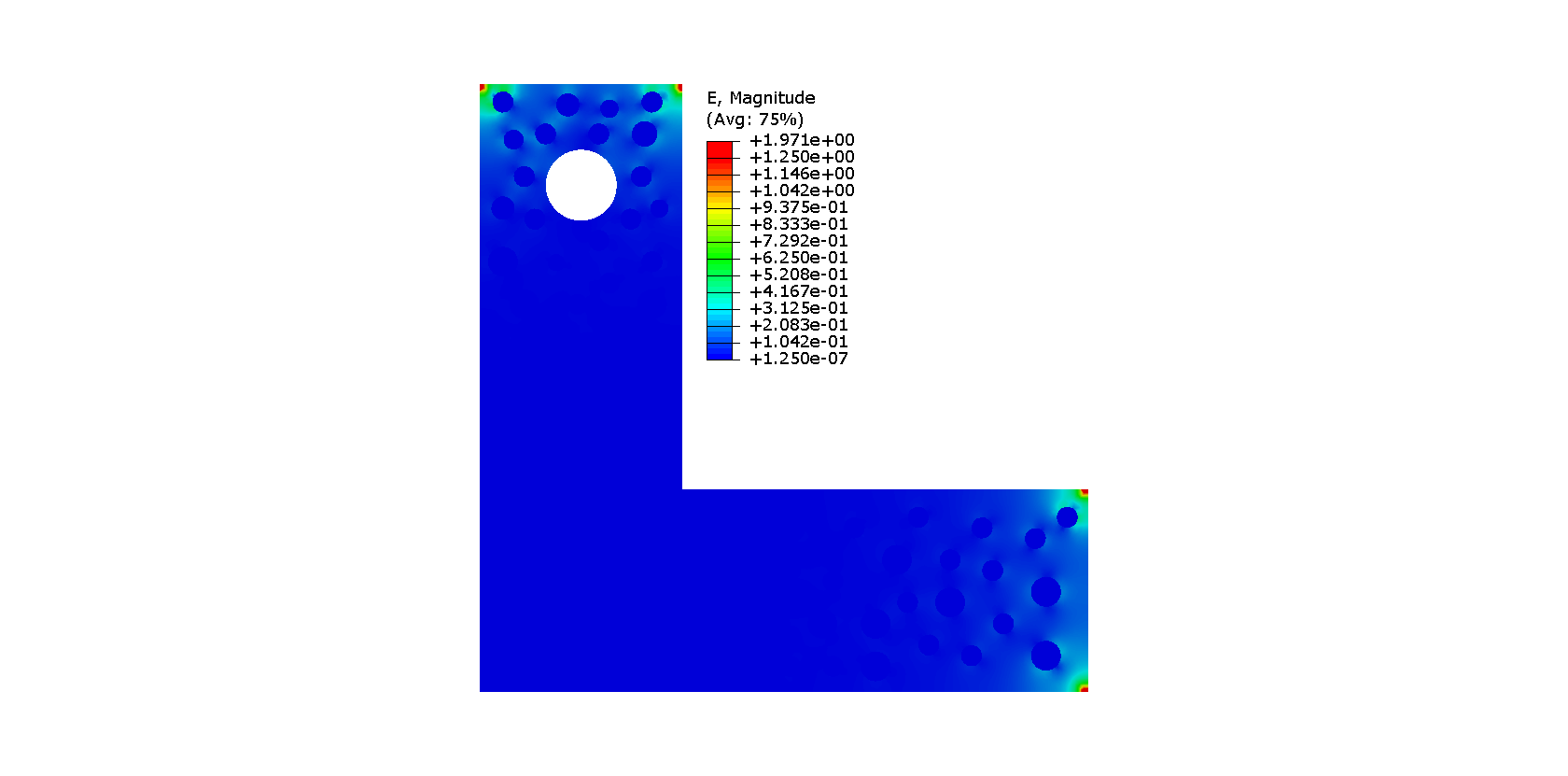}
      \caption{$\abs{\nabla \bar{u}(x)}$}
      \label{fig:geo4kklStrain}
   \end{subfigure}
\\
   \begin{subfigure}[hb]{0.39\textwidth}
      \includegraphics[trim={7in 1in 7in 1in},clip,width=\textwidth]{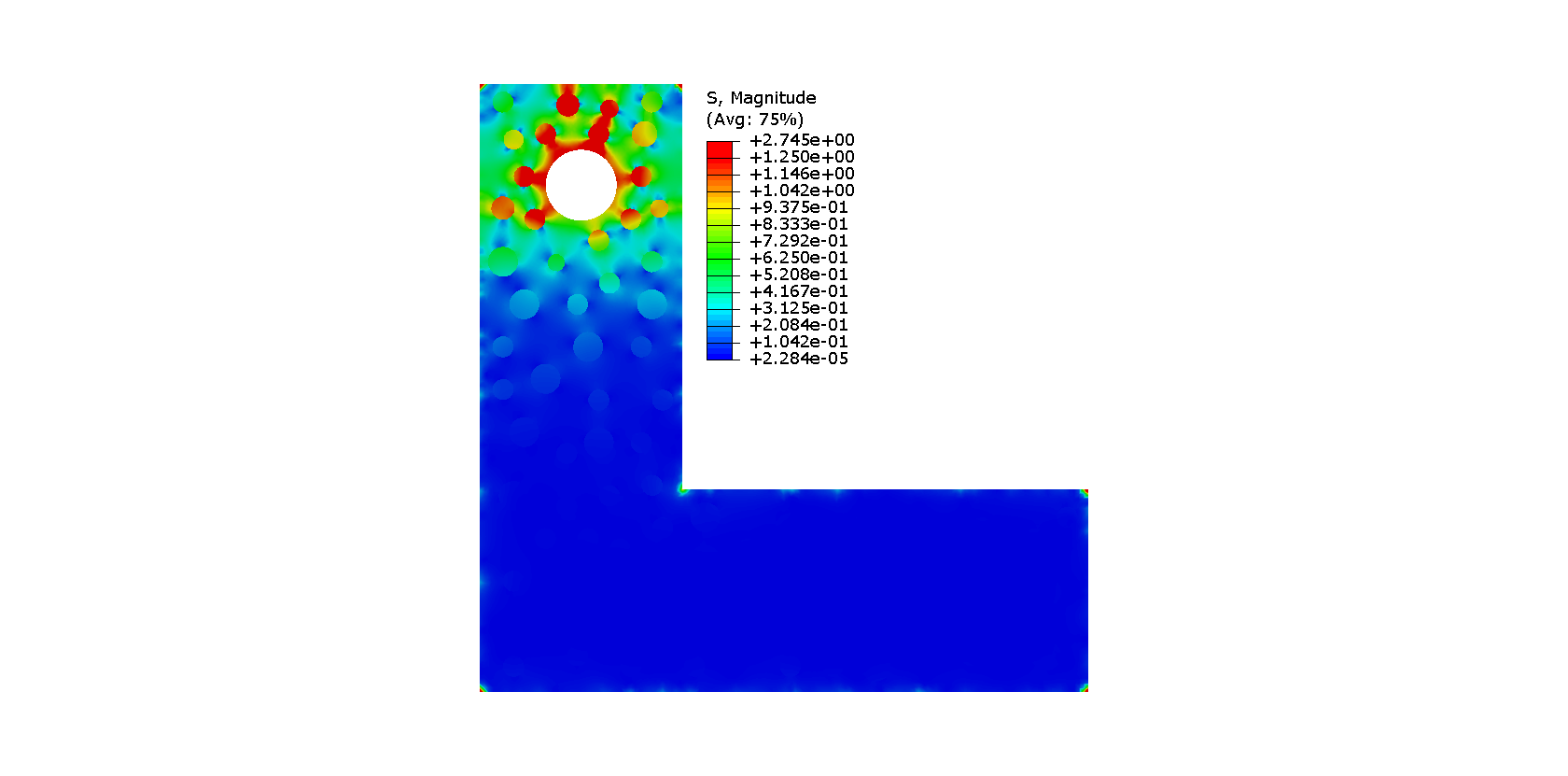}
      \caption{$\abs{c(x)\nabla \tilde{u}(x)}$}
      \label{fig:geo4eigStress}
   \end{subfigure}
   \begin{subfigure}[hb]{0.39\textwidth}
      \includegraphics[trim={7in 1in 7in 1in},clip,width=\textwidth]{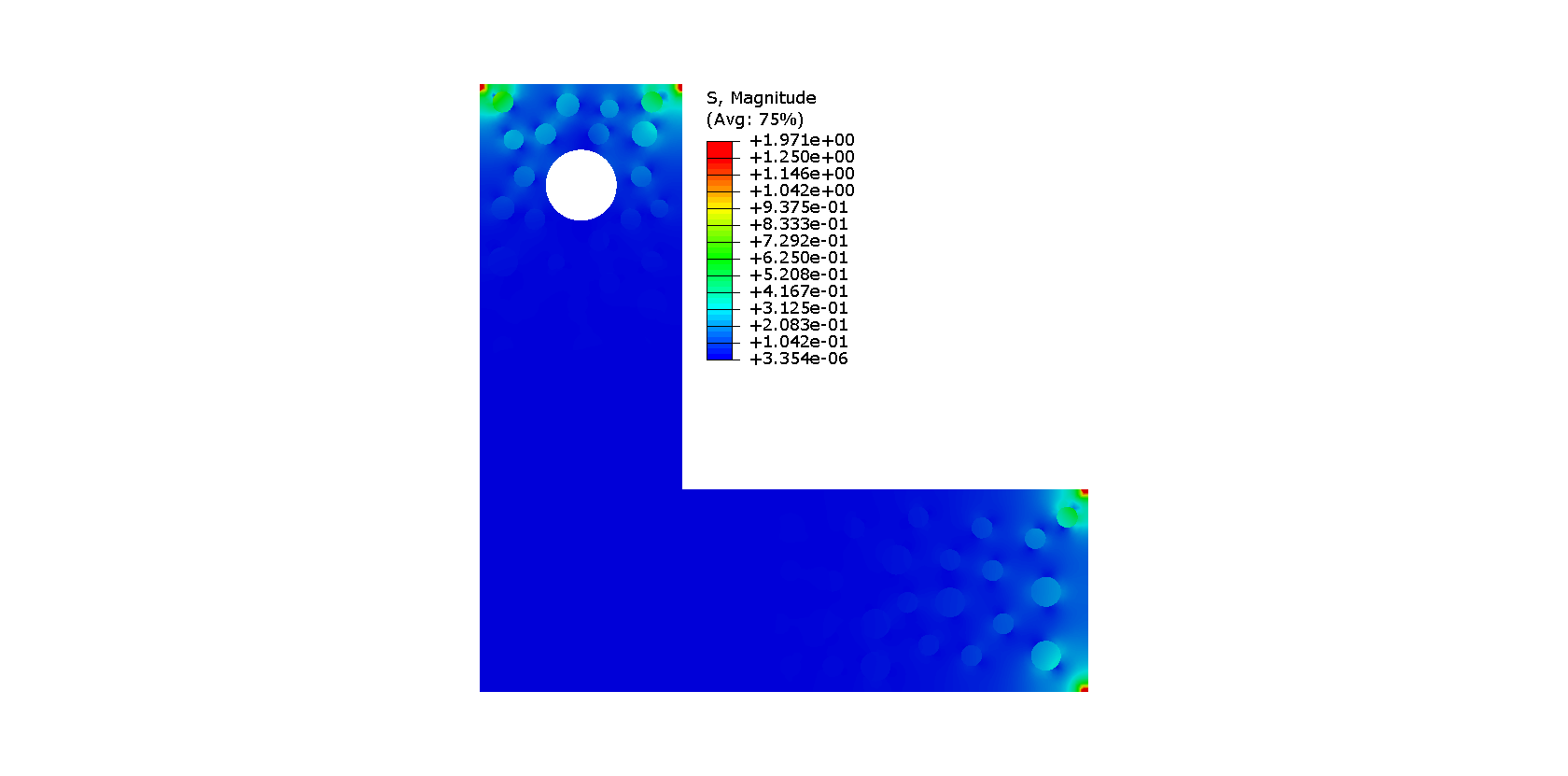}
      \caption{$\abs{c(x)\nabla \bar{u}(x)}$}
      \label{fig:geo4kklStress}
   \end{subfigure}
   \caption{Cross-section of Geometry 4.  The worst case solution (a) $\tilde{u}$ and the corresponding (c) strain and (e) stress compared to the ensemble averaged solution (b) $\bar{u}$ and the corresponding (d) strain and (f) stress.}
   \label{fig:geo4}
\end{figure}

Geometry 4 is an L shaped bracket with hole at the top leg as shown in Figure~\ref{fig:geo4Layout}.  The domain of interest $\omega$ is an annulus around the hole.  The bracket has long edges of length 6, short edges of length 4, and legs of width 2, the hole has radius 0.35, and $\omega$ radius 0.85.  The parameter $a$=12, and $b$=1.  The shear modulus $c(x)=1000$ in the inclusions and $c(x)=1$ in the matrix.  The ensemble average boundary displacement $\bar{g}(x)=\pm 0.1$ or $0$ as shown in Figure~\ref{fig:geo4Boundary}.  The computed fields are displayed in Figure~\ref{fig:geo4}.

\begin{table}[!ht]
%  \begin{center}
    \centering
    \begin{tabular}{| c | c | c | c |}
    \hline
    Geometry   &    $V$    &    \raisebox{0pt}[10pt][0pt]{$\bar{P}_N$}   &    $\hat{g}_n$, $\hat{g}_n'$ \\
    \hline
    1                & 0.604    &    0.136                                                           &    20, 17\\
    \hline
    2                & 0.695    &    0.149                                                           &     36, 25\\
    \hline
    3                 & 0.576    &   0.029                                                           &      65, 63\\
    \hline
    4                & 0.686     &   0.072                                                           &     54, 31 \\
    \hline
    \end{tabular}
%  \end{center}
  \caption{Computations of energy concentration for each geometry.  The final column shows the number of functions used from the closed form KKL expansion given by \eqref{eq:kklSolutions}}
  \label{tab:numerics}
\end{table}

%%%%%%%%%%%%%%%%%%%%%%%%%%%%%%%%%%%%%%%%%%%%%%%%%%%%
%%%%%%%%%%%%%%%%%%%%%%%%%%%%%%%%%%%%%%%%%%%%%%%%%%%%
%%%%%%%%%%%%%%%%%%%%%%%%%%%%%%%%%%%%%%%%%%%%%%%%%%%%

%%%%%%%%%%%%%%%%%%%%%%%%%%%%%%%%%%%%%%%%%%%%%%%%%%%%
%%%%%%%%%%%%%%%%%%%%%%%%%%%%%%%%%%%%%%%%%%%%%%%%%%%%
%%%%%%%%%%%%%%%%%%%%%%%%%%%%%%%%%%%%%%%%%%%%%%%%%%%%

\begin{figure}[H]
   \centering
   \begin{subfigure}[hb]{0.25\textwidth}
      \includegraphics[trim={7.5in 1.2in 10in 1.2in},clip,width=\textwidth]{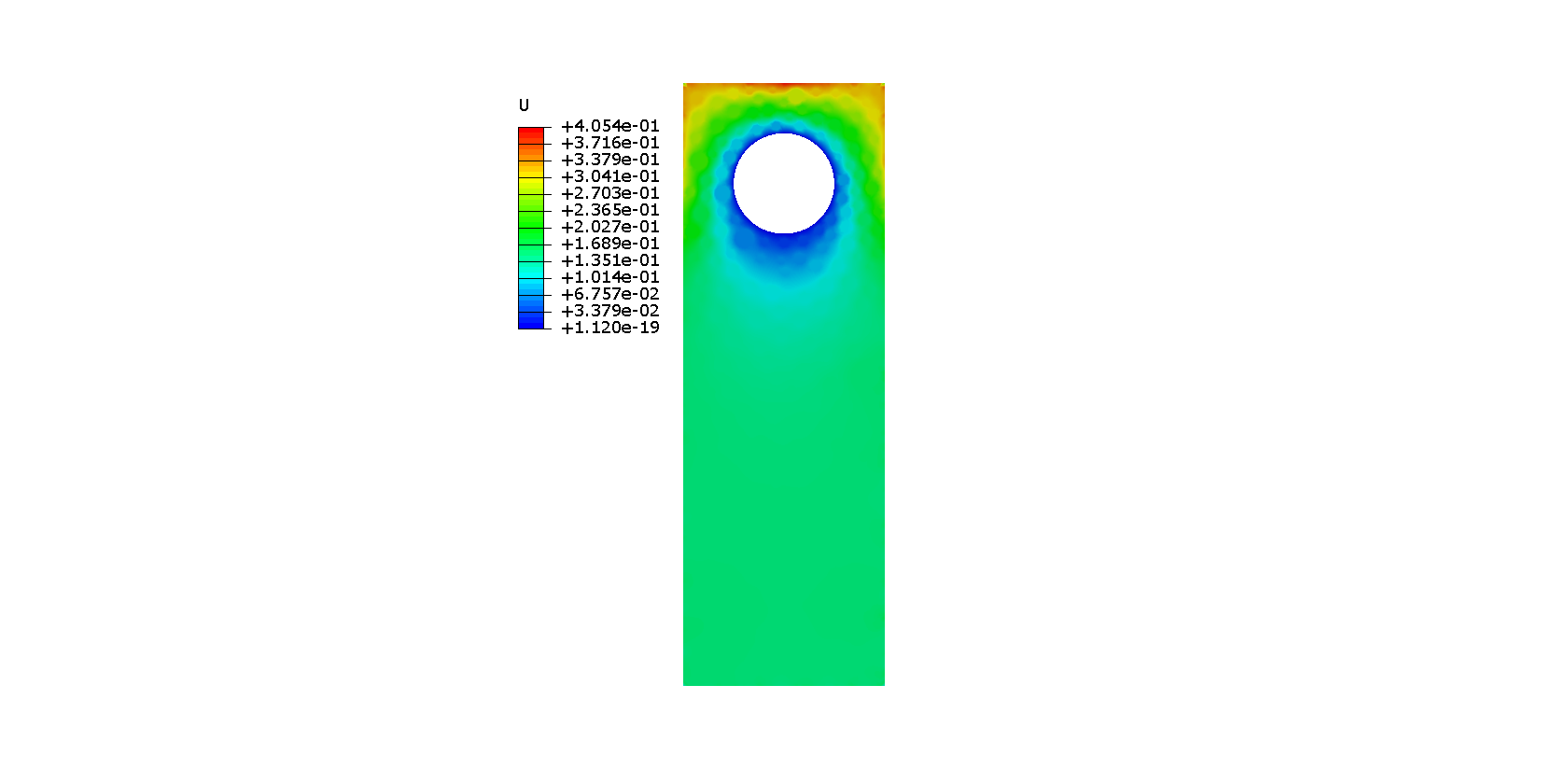}
      \caption{$\tilde{u}(x)$}
      \label{fig:geo2eig}
   \end{subfigure}
\qquad
   \begin{subfigure}[hb]{0.25\textwidth}
      \includegraphics[trim={7.5in 1.2in 10in 1.2in},clip,width=\textwidth]{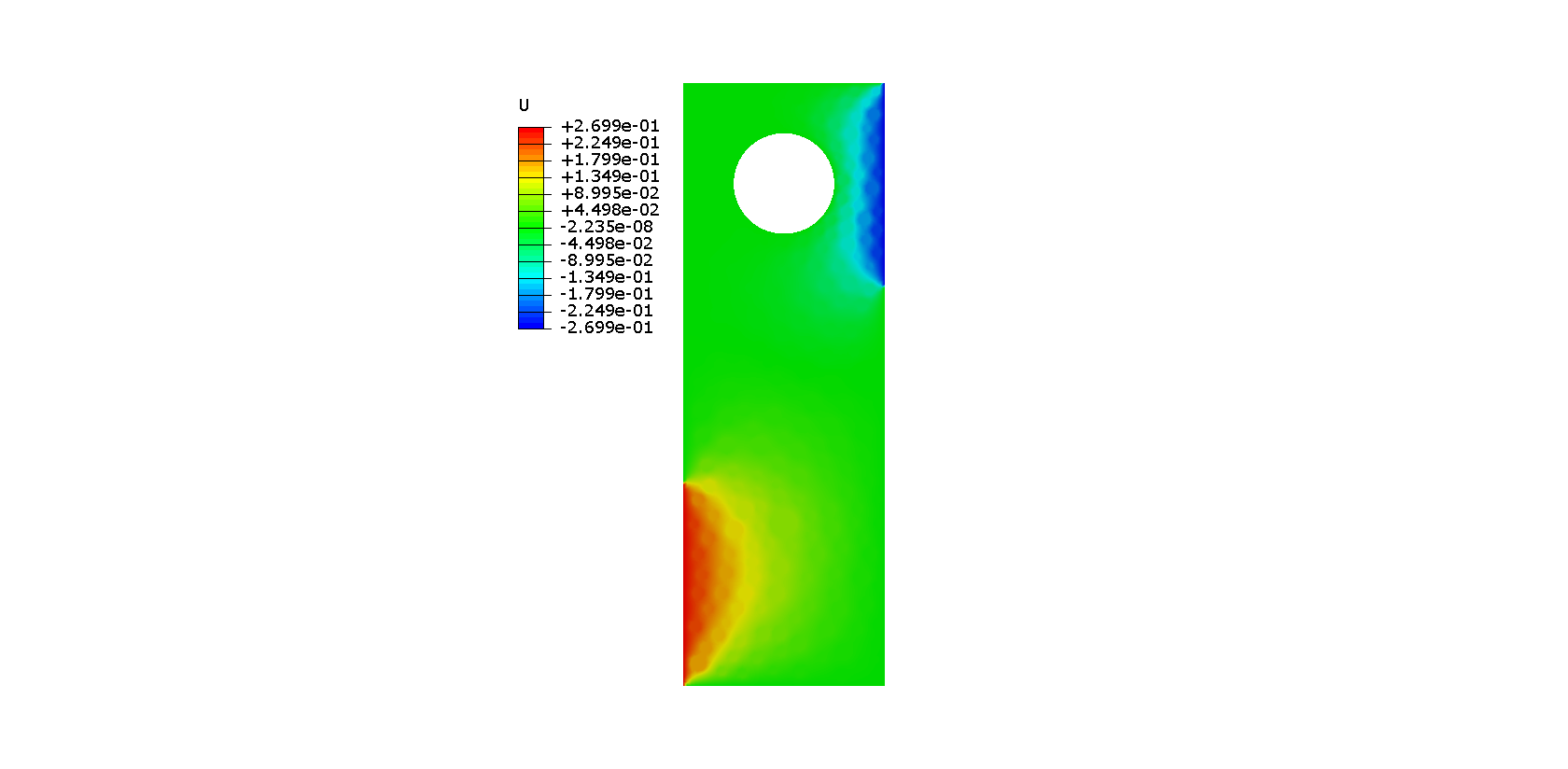}
      \caption{$\bar{u}(x)$}
      \label{fig:geo2kklAvg}
   \end{subfigure}
\\
   \begin{subfigure}[hb]{0.25\textwidth}
      \includegraphics[trim={7.5in 1.2in 10in 1.2in},clip,width=\textwidth]{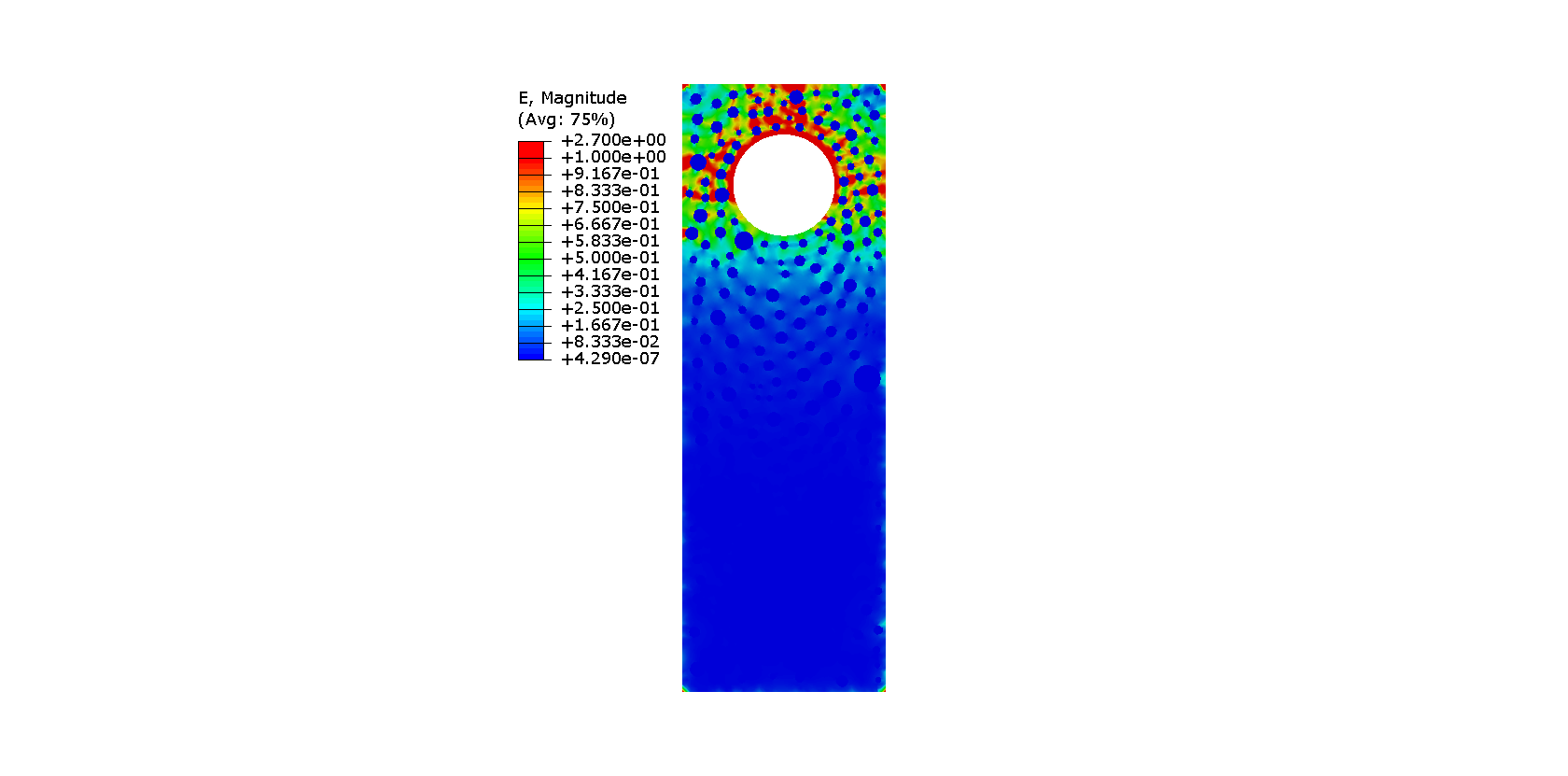}
      \caption{$\abs{\nabla \tilde{u}(x)}$}
      \label{fig:geo2eigStrain}
   \end{subfigure}
\qquad
   \begin{subfigure}[hb]{0.25\textwidth}
      \includegraphics[trim={7.5in 1.2in 10in 1.2in},clip,width=\textwidth]{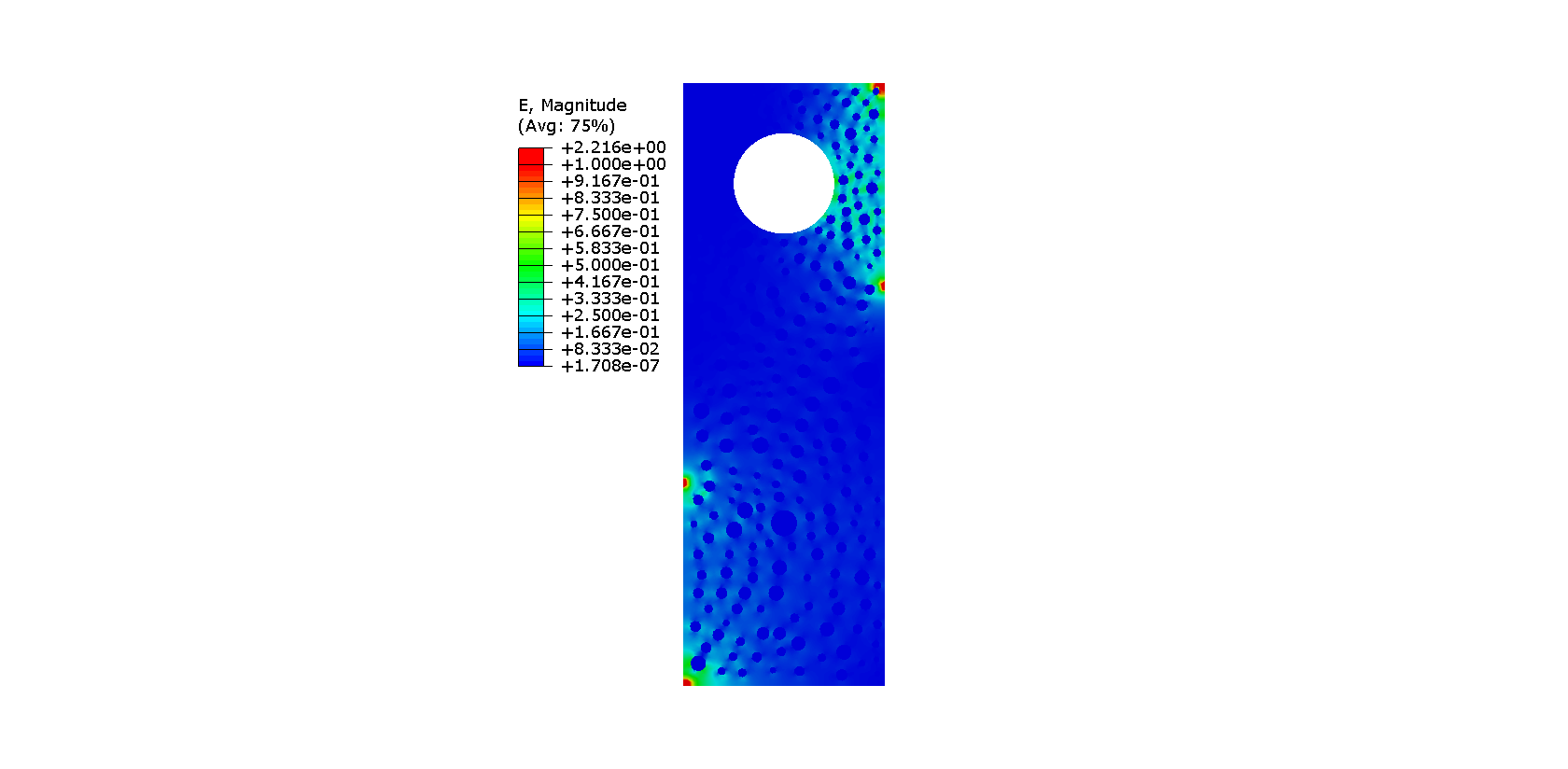}
      \caption{$\abs{\nabla \bar{u}(x)}$}
      \label{fig:geo2kklStrain}
   \end{subfigure}
\\
   \begin{subfigure}[hb]{0.25\textwidth}
      \includegraphics[trim={7.5in 1.2in 10in 1.2in},clip,width=\textwidth]{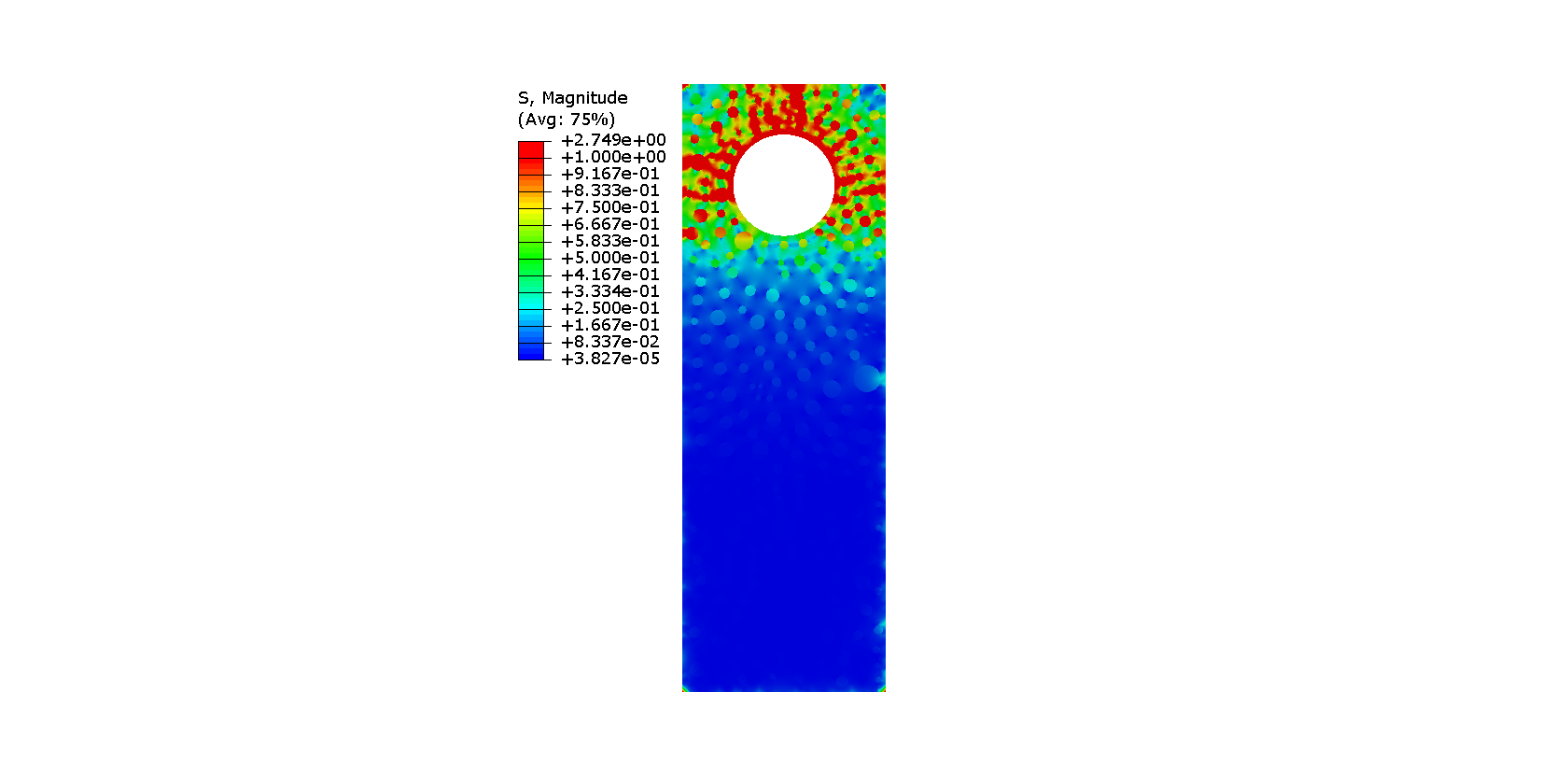}
      \caption{$\abs{c(x)\nabla \tilde{u}(x)}$}
      \label{fig:geo2eigStress}
   \end{subfigure}
\qquad
   \begin{subfigure}[hb]{0.25\textwidth}
      \includegraphics[trim={7.5in 1.2in 10in 1.2in},clip,width=\textwidth]{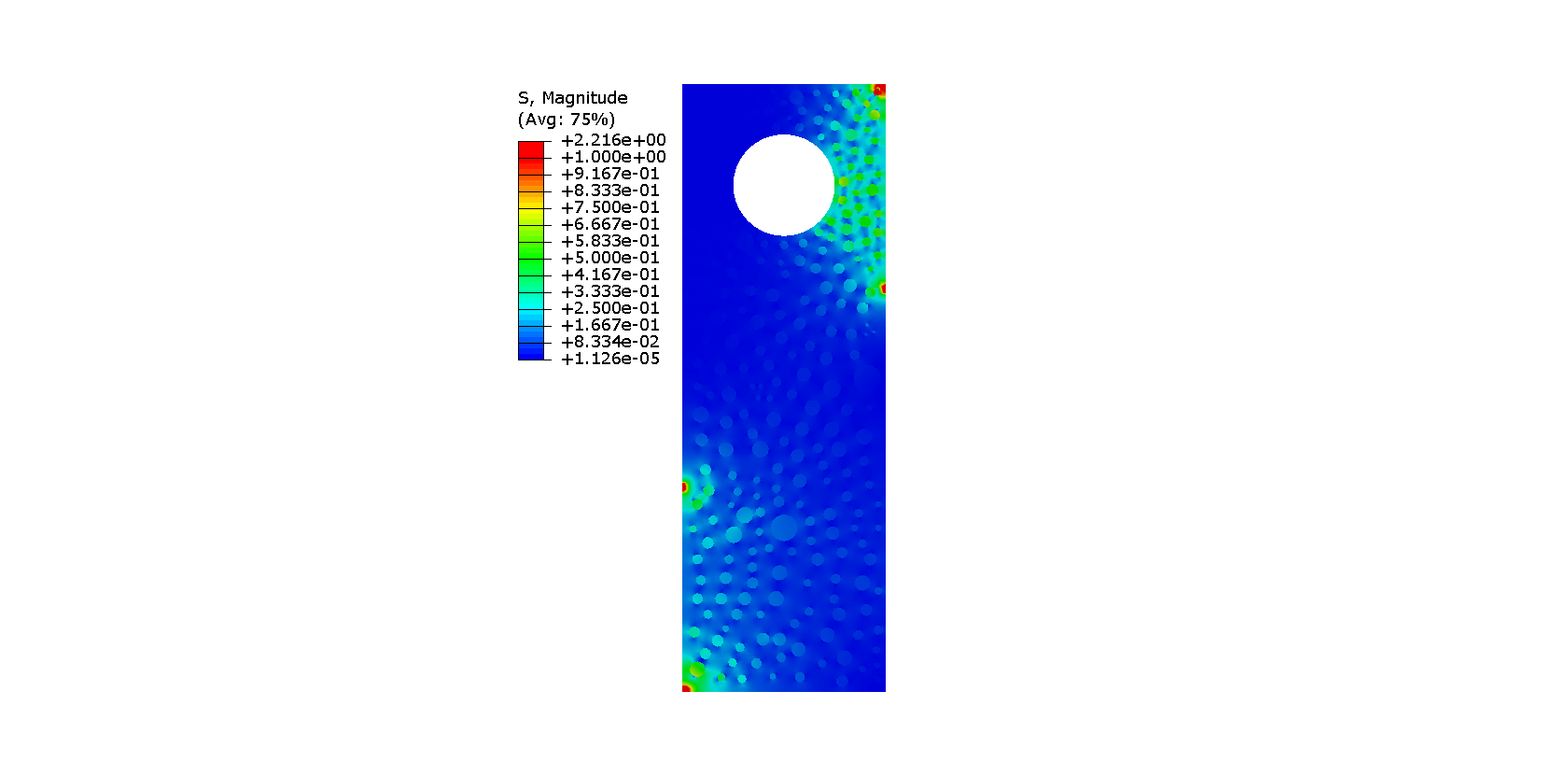}
      \caption{$\abs{c(x)\nabla \bar{u}(x)}$}
      \label{fig:geo2kklStress}
   \end{subfigure}
   \caption{Cross-section of Geometry 2.  The worst case solution (a) $\tilde{u}$ and the corresponding (c) strain and (e) stress compared to the ensemble averaged solution (b) $\bar{u}$ and the corresponding (d) strain and (f) stress.}
   \label{fig:geo2}
\end{figure}

%%%%%%%%%%%%%%%%%%%%%%%%%%%%%%%%%%%%%%%%%%%%%%%%%%%%
%%%%%%%%%%%%%%%%%%%%%%%%%%%%%%%%%%%%%%%%%%%%%%%%%%%%
%%%%%%%%%%%%%%%%%%%%%%%%%%%%%%%%%%%%%%%%%%%%%%%%%%%%

%%%%%%%%%%%%%%%%%%%%%%%%%%%%%%%%%%%%%%%%%%%%%%%%%%%%
%%%%%%%%%%%%%%%%%%%%%%%%%%%%%%%%%%%%%%%%%%%%%%%%%%%%
%%%%%%%%%%%%%%%%%%%%%%%%%%%%%%%%%%%%%%%%%%%%%%%%%%%%

\subsection{Observations}
Figures \ref{fig:geo2eig}, \ref{fig:geo3eig}, and \ref{fig:geo4eig} show, as expected, that the worst case load concentrates energy around the portion of the boundary near the domain of interest. The overall shape of the concentration for the worst case load is influenced by the placement of the heterogeneities.

The expected energy concentration  $\bar{P}_N$ is seen to be well below the worst case energy concentration  $V$, as seen in Table~\ref{tab:numerics}.  This large deviation is expected since the Markovian covariance leads to periodic boundary fluctuations \eqref{eq:kklSolutions} which distribute the random loads equally along the boundary.  A less well-behaved covariance could close this gap somewhat; however, the expected energy concentration falls considerably below the largest possible energy concentration associated with a particular ensemble of loads.

%%%%%%%%%%%%%%%%%%%%%%%%%%%%%%%%%%%%%%%%%%%%%%%%%%%%
%%%%%%%%%%%%%%%%%%%%%%%%%%%%%%%%%%%%%%%%%%%%%%%%%%%%
%%%%%%%%%%%%%%%%%%%%%%%%%%%%%%%%%%%%%%%%%%%%%%%%%%%%

\section{Conclusions}
We present a novel method for computing the worst case boundary load that imparts the maximum possible fraction of total energy onto a prescribed domain of interest contained within a composite structure.  The method works for both Neumann or Dirichlet boundary conditions.  The solution of the  maximal energy concentration problem is given by an eigenvalue eigenfunction pair for the concentration eigenvalue problem.  To illustrate the ideas this method is applied to four distinct geometries in the context of  anti-plane shear.  The the worst case loads are compared with the expected energy concentration of a random ensemble of loadings specified by a Markovian covariance with closed form solution using a KKL expansion.  The expected energy penetration of the random loadings for Markovian covariance were seen to be substantially less than the maximal energy concentration associated with the worst case load.   
The method presented here provides a novel way to recover a rigorous upper bound on how 'bad' a random ensemble of loadings could be.  The methods developed here can be used in evaluating and identifying worst case random loads with large energy concentration.  

Alternatively one can apply these methods to evaluate the ability of a composite geometry, such as a functionally graded composite, to prevent energy penetration to a particular interior subdomain. This approach can then be used as a design tool to find the design most resistant to the worst case load. 

\section{Acknowledgements}
This work supported in part through NSF Grant DMS-1211066 and by the Air Force Research Laboratory under University of Dayton Research Institute Contract FA8650-10-D-5011.

\end{document}